\documentclass[reqno]{amsart}
\usepackage{graphicx}
\usepackage{indentfirst,csquotes}
\usepackage{amssymb}
\usepackage{amsmath}
\usepackage{amsthm}
\usepackage{mathtools}
\usepackage{commath}
\usepackage{algorithm}
\usepackage{algpseudocode}
\usepackage{glossaries-extra}  
\usepackage{comment}
\usepackage{pgfplots}
\usepackage{array}
\usepackage{bigstrut}
\usepackage{graphicx} 
\usepackage{subcaption}
\usepackage{tikz}  
\usetikzlibrary{3d}
\usetikzlibrary{shapes.geometric}
\usepackage{hyperref}
\usepackage{fancyhdr}
\usepackage{xcolor}
\definecolor{linkblue}{rgb}{0,0.1,0.6}
\definecolor{citegreen}{rgb}{0,0.25,0.15}
\definecolor{linkred}{rgb}{0.8,0,0.005}
\definecolor{mailviolet}{rgb}{0.3,0,0.35}
\definecolor{tumblue}{rgb}{0,0.396,0.741}
\definecolor{tumorange}{RGB}{227,114,34}
\definecolor{tumgreen}{RGB}{162,173,0}
\definecolor{tumgrey}{RGB}{128,128,128}
\definecolor{tumLblue}{RGB}{100,160,200}
\definecolor{tumLgrey}{RGB}{218,215,203}
\definecolor{tumLgrey2}{RGB}{170,155,152}
\definecolor{tumLgrey3}{RGB}{204,204,204}

\begin{document}

\title{Parametric Model Order Reduction by Box Clustering with Applications in Mechatronic Systems}

\author{Juan Angelo Vargas-Fajardo$^{*,1}$, Diana Manvelyan-Stroot$^2$, Catharina Czech$^2$, Pietro Botazzoli$^3$, and Fabian Duddeck$^1$}

\date{\today \\[2mm]
    $^1$TUM School of Engineering and Design, Technical University of Munich, Arcisstr. 21, 80333 Munich\\%
    $^2$Siemens AG, Garching\\
    $^3$Siemens AG, Erlangen\\
    $^*$Corresponding author (juan.vargas@tum.de)\\
}

\begin{abstract}
High temperatures and structural deformations can compromise the functionality and reliability of new components for mechatronic systems. Therefore, high-fidelity simulations (HFS) are employed during the design process, as they enable a detailed analysis of the thermal and structural behavior of the system. However, such simulations are both computationally expensive and tedious, particularly during iterative optimization procedures. Establishing a parametric reduced order model (pROM) can accelerate the design's optimization if the model can accurately predict the behavior over a wide range of material and geometric properties. However, many existing methods exhibit limitations when applied to wide design ranges.
\\
In this work, we introduce the parametric Box Reduction (pBR) method, a matrix interpolation technique that minimizes the non-physical influence of training points due to the large parameter ranges. For this purpose, we define a new interpolation function that computes a local weight for each design variable and integrates them into the global function. Furthermore, we develop an intuitive clustering technique to select the training points for the model, avoiding numerical artifacts from distant points. Additionally, these two strategies do not require normalizing the parameter space and handle every property equally. The effectiveness of the pBR method is validated through two physical applications: structural deformation of a cantilever Timoshenko beam and heat transfer of a power module of a power converter. The results demonstrate that the pBR approach can accurately capture the behavior of mechatronic components across large parameter ranges without sacrificing computational efficiency.

\end{abstract}

\pagestyle{fancy}
\lhead{}
\rhead{}
  \fancyhead[OL]{\small{J.A. Vargas-Fajardo, D. Manvelyan-Stroot, C. Czech, P. Botazzoli and F. Duddeck}}

  \fancyhead[EL]{\small{Parametric Model Order Reduction by Box Clustering with Applications in Mechatronic Systems}}
  \rhead{\thepage}
 \cfoot{}

\maketitle
\let\thefootnote\relax

\section{Introduction}

Mechanical and thermal simulations play a role in most engineering industries, such as energy, aerospace, and automotive. They are crucial to the design and thermal management of mechatronic systems such as power modules, which are often prone to damaging heat distributions. High temperatures can lead to uncontrolled and non-recoverable structural deformations, eventually compromising the functionality and reliability of the device. Therefore, during operation, a cooling system maintains the temperature within the optimal range to avoid damage. The specific choice of the design parameters, such as material and geometric properties, can strongly affect the temperature distribution of the semiconductor chips in a power module. Numerical simulation allows a better understanding of this phenomenon during the development phase. 

Such high-fidelity simulations (HFS), typically based on the finite element (FEM) or finite volume method (FVM), can accurately evaluate the relevant design quantities of detailed three-dimensional virtual representations.

However, despite the constant increase in computing resources, computationally expensive HFS hinders design engineers from quickly evaluating design adaptations. This is particularly true during design optimization, where multiple analyses with varying design variables are executed to find the best performance. Therefore, parametric reduced order models (pROM) have been developed to replace HFS and speed up these iterative procedures. In recent years, pROMs have been extensively researched, leading to the development of new algorithms for generating surrogate models that can perform fast simulations \cite{Benner2021_V1,Baur2014}. Here, we focus on linear pROMs as the decisive structural mechanics and heat transfer equations for the application case are of a linear nature. Moreover, an intrusive reduction technique is chosen, employing the system matrices of the HFS. These approaches often achieve low approximation errors, as underlying system knowledge is included in constructing the pROM. 

Krylov subspace methods \cite{Freund2003,salimbahrami2006order} are a common intrusive reduction technique. They approximate the system's transfer function by matching the moments derived from the bases of a Krylov subspace, see \cite{Feldmann1995, Freund1999, Liesen2012}. Since its development, Krylov subspace methods have been widely employed in structural mechanics, circuit simulation, and thermal problems, with first results documented in \cite{Häggblad1993, Freund2000, Salleras2005}, thereby highlighting their broad applicability. One of the earliest and most detailed research studies for electro-thermal simulations was conducted in \cite{Bechtold2005}, where a ROM with the Arnoldi method is generated. However, in \cite{Bechtold2005}, the underlying Krylov subspace is constructed based on constant system matrices and does not include any parameter dependencies.  

For parameter-dependent problems, Krylov subspace methods can be combined with matrix interpolation techniques to create parametric surrogate models \cite{Panzer2010}, consisting of three main stages: reducing full-order models (FOMs) into local ROMs using a Krylov subspace method, finding a global system representation for these ROMs by changing their basis, and finally, interpolating the reduced system matrices to evaluate new parameter points. Given the parametric nature of the model, the design parameter values need to be normalized before any of the previous steps are performed. Many normalization techniques exist, such as Min-Max, Z-Score, and Decimal scaling. A clear drawback is that the model highly depends on the normalization technique, especially for a high-dimensional parameter space \cite{He2021,Muhammad2022,Pandey2017}. 

Following its introduction in \cite{Panzer2010}, numerous studies have been conducted to improve its efficiency and applicability. In \cite{Hassan2023, Hassan2024}, the interpolation step is modified by using a Lagrange interpolation technique for power module applications. In \cite{Geuss2014,Geuss2016,Brunsch2017}, different matrix interpolation approaches for high-dimensional parameter spaces are introduced, which employ sparse grids. Another option for parametric models with Krylov methods is using a multivariate expansion to approximate the parametric transfer function \cite{Baur2011, Yuan2021}.

In contrast to intrusive pROMs, non-intrusive methods treat the HFS as a black-box, see \cite{Drofenik2007, Zhang2021, Kulkarni2023} for reduction techniques in similar application cases. Recently, machine learning algorithms, such as convolutional neural networks (CNN), have been used to approximate heat conduction problems \cite{Tadeparti2022, Peng2021, Rezasefat2023}. For instance, \cite{Stipsitz2022} extends this approach to three-dimensional thermal surrogate models for electronic systems. However, most studies focus on steady-state problems. 

This work presents two new methods for modifying the framework presented in \cite{Panzer2010}. The novelty of the proposed approach lies in its ability to generate pMORs without prior normalization of the parameter space, thereby eliminating the dependence on normalization functions and their associated hyperparameters. Furthermore, the generated model is accurate across wide parameter ranges, as it mitigates non-physical influences from training points without compromising computational efficiency. The methods are divided into three main stages: (i) system reduction, where with a set of training points reduced models are created with a Krylov Subspace Method; (ii) basis change, where all reduced models are projected into a global reduced space; (iii) matrix interpolation, where the new reduced matrices for any given combination of input parameters are computed. With this reduced matrices, fast simulations can be computed.

The first approach, parametric box interpolation, computes the matrix interpolation using a newly introduced weighting function that employs a subset of training points. The selection process is carried out using the box clustering algorithm, ensuring a consistent sampling without normalization. The second and primary method of this contribution, parametric box reduction, extends the box clustering to the basis change stage from \cite{Panzer2010}. It removes the nonphysical numerical influence from distant training points during the pROMs generation.

The structure of this work is as follows: Section \ref{ss:formulation} introduces the general problem formulation. Afterward, Section \ref{ss:pBR} presents the theoretical background of the new technique and provides a detailed explanation of the new methods. To assess the technique's performance, reduced models are generated for a structural problem of a 3D Timoshenko beam to illustrate the principles of the approach, and for the thermal behavior of a power module to show its potential in thermo-mechanical problems. Their benefits and drawbacks are discussed in terms of accuracy and computational time in Section \ref{ss:results}. Finally, Section \ref{ss:conclusion} summarizes the main findings and gives an outlook for future research paths. 
\section{Problem formulation}\label{ss:formulation}

For many application cases, PDEs can describe the most important physical phenomena of the problem at hand. The FEM and FVM are well-known, time-tested methods used to perform HFS to predict the behavior of the relevant quantities of the design. Both FEM and FVM discretize the domain and establish the characteristic system matrices. They create a linear time-invariant (LTI) system, which is often a first-order system, as shown in the following equation. 

\begin{equation}\label{eq:FOM_1st}
\begin{aligned}
    \mathbf{E}\left(\mu_i\right)\dot{u}\left(t,\mu_i\right) &= \mathbf{A}\left(\mu_i\right)u\left(t,\mu_i\right) + f\left(t,\mu_i\right), \\
    y\left(t,\mu_i\right) &= \mathbf{D}u\left(t,\mu_i\right) \, ,
\end{aligned}
\end{equation}

where $\mathbf{E}(\mu_i) \in \mathbb{R}^{n \times n}$ and $\mathbf{A}(\mu_i) \in \mathbb{R}^{n \times n}$ are the given system matrices, $f(\cdot,\mu_i) \in \mathbb{R}^{n}$ is the load vector, $\mathbf{D}(\mu_i) \in \mathbb{R}^{m \times n}$ the output matrix, with $n$ the number of degrees of freedom (i.e., dimension of the HFS) and $m$ the number of outputs. The parameter-dependent vectors $u(\cdot, \mu_i)$ and $y(\cdot, \mu_i)$ are the state and output variables, respectively. The system matrices are a function of $\mu_i$, which is a design point of the parameter set $\mathbb{P} = \{\mu_1, \ldots, \mu_k\}$. For better readability, from this point onward, a matrix or vector with subindex $i$ denotes the matrix evaluated at the parameter point $\mu_i$. In other cases, the discretization leads to an undamped second-order system

\begin{equation}\label{eq:FOM_2nd}
\begin{aligned}
    \mathbf{M}_i\ddot{u}_i\left(t\right) &= \mathbf{K}_iu_i\left(t\right) + f_i\left(t\right), \\
    y_i\left(t\right) &= \mathbf{D}u_i\left(t\right)
\end{aligned}
\end{equation}

where $\mathbf{M}_i$ and $\mathbf{K}_i$ are the mass and stiffness matrices respectively.
From the systems of Eq. \eqref{eq:FOM_1st}-\eqref{eq:FOM_2nd} and a numerical time integration scheme, the primary variable $u_i(t)$ at each time step is computed for the design point $\mu_i$. The design points have $p$ design variables $x_j$.

\begin{equation}\label{eq:ParameterSpace}
\begin{aligned}
    \mu_i = \left[x_1, x_2 , \dots x_p\right]^T \in \mathbb{P} \subset  \mathbb{R}^{p}\,.
\end{aligned}
\end{equation}

To train the model, the training set $\mathbb{T}$ is constructed by sampling $k$ training points $\mu_i \in \mathbb{P}$; $k$ system matrices of the $k$ corresponding HFMs are extracted, as we are employing an intrusive reduction technique. The pROM predicts the primary variable $u_{r,j}(t)$ for a new parameter point $\mu_j \in \mathbb{P} \setminus \mathbb{T}$, where the subindex $r$ indicates that it results from a pROM. 
\section{Parametric Box Reduction}\label{ss:pBR}

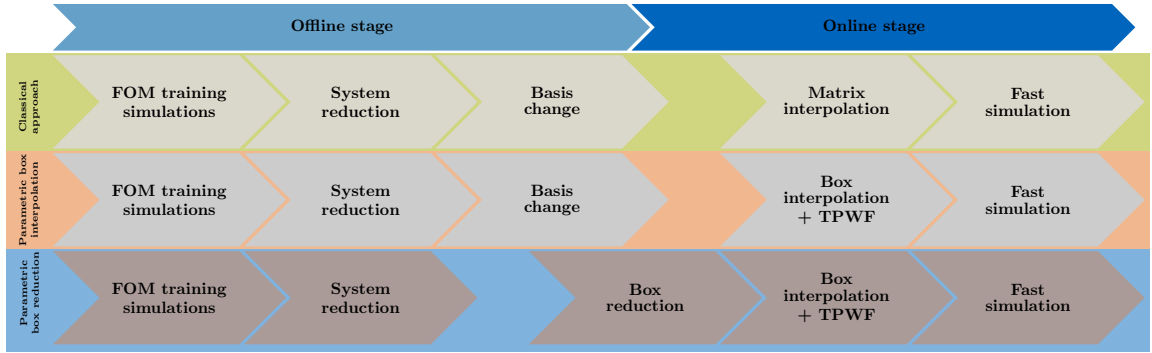
\begin{figure} [ht] 
    \centering
    \resizebox{0.95\textwidth}{!}{
    \begin{tikzpicture}
        \tikzstyle{mypolygon} = [draw = none, fill=tumLgrey]
        \tikzstyle{mypolygon3} = [draw = none, fill=tumgrey]
        \tikzstyle{mypolygon4} = [draw = none, fill=tumLblue]
        \tikzstyle{mypolygon8} = [draw = none, fill=tumgreen,opacity=0.5]
        \tikzstyle{mypolygon5} = [draw = none, fill=tumorange, opacity=0.5]
        \tikzstyle{mypolygon7} = [draw = none, fill=tumblue, opacity=0.5]
        \tikzstyle{mypolygon6} = [draw = none, fill=tumblue]
        \tikzstyle{mypolygon2} = [draw = none, fill=tumLgrey2]        
        \tikzstyle{mypolygon9} = [draw = none, fill=tumLgrey3]


        \draw[mypolygon5] (1, 5.95) -- (25.5,5.95) -- (25.5,8.05) -- (1,8.05)  -- cycle;
        \node[rotate=90] at (1.5, 7) {\shortstack{\small \tiny\textbf{Parametric box} \\ \tiny \textbf{interpolation}}}; 

        \draw[mypolygon8] (1, 8.05) -- (25.5,8.05) -- (25.5,10.15) -- (1,10.15)  -- cycle;
        \node[rotate=90] at (1.5, 9) {\shortstack{\tiny\textbf{Classical} \\ \tiny\textbf{approach}\cite{Panzer2010}}};

        \draw[mypolygon7] (1, 3.75) -- (25.5, 3.75) -- (25.5,5.95) -- (1,5.95)  -- cycle;
        \node[rotate=90] at (1.5, 5) {\shortstack{\tiny\textbf{Parametric} \\ \tiny\textbf{box reduction}}};

        \def\x{2}
        \def\y{8.1}
        \draw[mypolygon] (\x, \y) -- (\x + 4, \y) -- (\x + 5, \y + 1) -- (\x + 4, \y + 2) -- (\x , \y + 2) -- (\x + 1, \y + 1) -- cycle;
        \node at (\x + 2.5, \y + 1) {\shortstack{\textbf{FOM training} \\ \textbf{simulations}}};

        \def\x{6.1}
        \def\y{8.1}
        \draw[mypolygon] (\x, \y) -- (\x + 4, \y) -- (\x + 5, \y + 1) -- (\x + 4, \y + 2) -- (\x , \y + 2) -- (\x + 1, \y + 1) -- cycle;
        \node at (\x + 2.5, \y + 1) {\shortstack{\textbf{System} \\ \textbf{reduction}}};

        \def\x{10.2}
        \def\y{8.1}
        \draw[mypolygon] (\x, \y) -- (\x + 4, \y) -- (\x + 5, \y + 1) -- (\x + 4, \y + 2) -- (\x , \y + 2) -- (\x + 1, \y + 1) -- cycle;
        \node at (\x + 2.5, \y + 1) {\shortstack{\textbf{Basis} \\ \textbf{change}}};

        \def\x{16.3}
        \def\y{8.1}
        \draw[mypolygon] (\x, \y) -- (\x + 4, \y) -- (\x + 5, \y + 1) -- (\x + 4, \y + 2) -- (\x , \y + 2) -- (\x + 1, \y + 1) -- cycle;
        \node at (\x + 2.5, \y + 1) {\shortstack{\textbf{Matrix} \\ \textbf{interpolation}}};

        \def\x{20.4}
        \def\y{8.1}
        \draw[mypolygon] (\x, \y) -- (\x + 4, \y) -- (\x + 5, \y + 1) -- (\x + 4, \y + 2) -- (\x , \y + 2) -- (\x + 1, \y + 1) -- cycle;
        \node at (\x + 2.5, \y + 1) {\shortstack{\textbf{Fast} \\ \textbf{simulation}}};


        \def\x{2}
        \def\y{6}
        \draw[mypolygon9] (\x, \y) -- (\x + 4, \y) -- (\x + 5, \y + 1) -- (\x + 4, \y + 2) -- (\x , \y + 2) -- (\x + 1, \y + 1) -- cycle;
        \node at (\x + 2.5, \y + 1) {\shortstack{\textbf{FOM training} \\ \textbf{simulations}}};

        \def\x{6.1}
        \def\y{6}
        \draw[mypolygon9] (\x, \y) -- (\x + 4, \y) -- (\x + 5, \y + 1) -- (\x + 4, \y + 2) -- (\x , \y + 2) -- (\x + 1, \y + 1) -- cycle;
        \node at (\x + 2.5, \y + 1) {\shortstack{\textbf{System} \\ \textbf{reduction}}};

        \def\x{10.2}
        \def\y{6}
        \draw[mypolygon9] (\x, \y) -- (\x + 4, \y) -- (\x + 5, \y + 1) -- (\x + 4, \y + 2) -- (\x , \y + 2) -- (\x + 1, \y + 1) -- cycle;
        \node at (\x + 2.5, \y + 1) {\shortstack{\textbf{Basis} \\ \textbf{change}}};

        \def\x{16.3}
        \def\y{6}
        \draw[mypolygon9] (\x, \y) -- (\x + 4, \y) -- (\x + 5, \y + 1) -- (\x + 4, \y + 2) -- (\x , \y + 2) -- (\x + 1, \y + 1) -- cycle;
        \node at (\x + 2.5, \y + 1) {\shortstack{\textbf{Box} \\ \textbf{interpolation} \\\textbf{+ TPWF} }};

        \def\x{20.4}
        \def\y{6}
        \draw[mypolygon9] (\x, \y) -- (\x + 4, \y) -- (\x + 5, \y + 1) -- (\x + 4, \y + 2) -- (\x , \y + 2) -- (\x + 1, \y + 1) -- cycle;
        \node at (\x + 2.5, \y + 1) {\shortstack{\textbf{Fast} \\ \textbf{simulation}}};


        \def\x{2}
        \def\y{3.9}
        \draw[mypolygon2] (\x, \y) -- (\x + 4, \y) -- (\x + 5, \y + 1) -- (\x + 4, \y + 2) -- (\x , \y + 2) -- (\x + 1, \y + 1) -- cycle;
        \node at (\x + 2.5, \y + 1) {\shortstack{\textbf{FOM training} \\ \textbf{simulations}}};

        \def\x{6.1}
        \def\y{3.9}
        \draw[mypolygon2] (\x, \y) -- (\x + 4, \y) -- (\x + 5, \y + 1) -- (\x + 4, \y + 2) -- (\x , \y + 2) -- (\x + 1, \y + 1) -- cycle;
        \node at (\x + 2.5, \y + 1) {\shortstack{\textbf{System} \\ \textbf{reduction}}};

        \def\x{12.2}
        \def\y{3.9}
        \draw[mypolygon2] (\x, \y) -- (\x + 4, \y) -- (\x + 5, \y + 1) -- (\x + 4, \y + 2) -- (\x , \y + 2) -- (\x + 1, \y + 1) -- cycle;
        \node at (\x + 2.5, \y + 1) {\shortstack{\textbf{Box} \\ \textbf{reduction}}};

        \def\x{16.3}
        \def\y{3.9}
        \draw[mypolygon2] (\x, \y) -- (\x + 4, \y) -- (\x + 5, \y + 1) -- (\x + 4, \y + 2) -- (\x , \y + 2) -- (\x + 1, \y + 1) -- cycle;
        \node at (\x + 2.5, \y + 1) {\shortstack{\textbf{Box} \\ \textbf{interpolation} \\\textbf{+ TPWF} }};

        \def\x{20.4}
        \def\y{3.9}
        \draw[mypolygon2] (\x, \y) -- (\x + 4, \y) -- (\x + 5, \y + 1) -- (\x + 4, \y + 2) -- (\x , \y + 2) -- (\x + 1, \y + 1) -- cycle;
        \node at (\x + 2.5, \y + 1) {\shortstack{\textbf{Fast} \\ \textbf{simulation}}};

        \def\x{2}
        \def\y{10.2}
        \draw[mypolygon4] (\x, \y) -- (\x + 12.3, \y) -- (\x + 12.8, \y + 0.5) -- (\x + 12.3, \y + 1) -- (\x , \y + 1) -- (\x + .5, \y + 0.5) -- cycle;
        \node at (\x + 6.2, \y + 0.5) {\textbf{Offline stage}};

        \def\x{14.4}
        \def\y{10.2}
        \draw[mypolygon6] (\x, \y) -- (\x + 10.3, \y) -- (\x + 10.8, \y + 0.5) -- (\x + 10.3, \y + 1) -- (\x , \y + 1) -- (\x + .5, \y + 0.5) -- cycle;
        \node at (\x + 5.2, \y + 0.5) {\textbf{Online stage}};

    \end{tikzpicture}
    }
    \caption{Schematic representation of the main reduction steps of the classical approach \cite{Panzer2010} (green), parametric box interpolation (orange), and parametric box reduction (blue). pBI only changes the matrix interpolation step, while pBR modifies the basis change step as well.}\label{fig:pBR2}  
\end{figure}

This section describes the novel procedure for generating parametric reduced order models from high-fidelity simulations. Initially, a technique combining Krylov subspace methods with matrix interpolation to create pROM is introduced, referred to as the classical approach in this work \cite{Panzer2010}. This method has three main steps: system reduction, basis change, and matrix interpolation. The first two are carried out during the offline or training phase, while the latter occurs in an online phase. This method generates surrogate models, with which the designer can perform multiple fast simulations for different designs without creating additional HFMs. A process schematic is depicted in \textbf{green} in Fig. \ref{fig:pBR2}; this technique serves as the reference framework for the novel method. Like many other pMOR techniques, the classical approach encounters challenges in accurately predicting behavior when dealing with large parameter ranges. These complications are especially true when materials have not been selected during the initial design phase because the properties between two different materials can vary significantly.

Given the nature of the parameter space (i.e., large parameter ranges), the values of the design variables values need to be normalized before any calculation is performed. This step is essential to treat every parameter equally, thus reducing any non-physical influence on the results. Many normalization techniques exist, such as Min-Max, Z-Score, and Decimal scaling. A clear drawback is that the model highly depends on the normalization technique, especially for a high-dimensional parameter space \cite{He2021,Muhammad2022,Pandey2017}. Additionally, many of these techniques require tuning hyperparameters, leading to a higher effort during the model creation. Similarly, many interpolation weight functions $\omega$ have tunable parameters, which affect the model's accuracy. 

Considering these issues, we propose the \textit{Tensor Product Weight Function} (TPWF) (Sect. \ref{ss:TPWF}) as an alternative to the classical approaches. Normalizing the parameter space is not required to compute the weight function. Furthermore, TPWF treats each parameter equally, avoiding non-physical influence in the pROM. TPWF computes a local weight function for each dimension, and similarly to high-dimension shape functions, they are combined to have a global weight function. Although related ideas are presented in \cite{Geuss2014, Hassan2023, Hassan2024}, the proposed approach has not been used for computing matrix interpolation's weight functions.

Moreover, another issue arising from large parameter ranges is the dissimilarity of physical behavior between the training points, which reduces the accuracy of the prediction. A common practice to improve the model's accuracy is adding more training points. Although more information is added, poor model performance may occur due to the significant difference between the training points behavior. To avoid this issue, clustering techniques have been developed to select adequate interpolation points. The authors in \cite{Schopper2024} introduce different approaches to partition the parameter space in regions where all local reduced bases are consistent with specific predefined metrics. In \cite{Yuan2021b}, the authors chose only the nearby ROMs for the interpolation. The selected training points are identified with a $k$-Nearest Neighbor algorithm (kNN).

Nevertheless, the parameter space range must be normalized before the clustering. Consequently, this leads to the same issues mentioned before. Considering this, we present \textit{Box Clustering} (Sect. \ref{ss:BC}): a clustering technique to choose the points that compute the reduced system matrices of the new parameter point. It finds all the training points that create a box surrounding one parameter. This scheme does not require a normalized parameter space. Even if the ranges are normalized, the selected points remain unchanged. Furthermore, it ensures that no design variable has any additional artificial influence on the model due to the clustering technique. 

By combining these two methods, we introduce the \textit{Parametric Box Interpolation} (pBI) (Sect. \ref{ss:BI}) technique. It employs TPWF and box clustering during the matrix interpolation step for computing the new reduced system, while the other steps remain the same as in the classical approach. Fig. \ref{fig:pBR2} illustrates in \textbf{orange} the main steps of the new procedure. 

Although pBI reduces the non-physical influence of training points in the online stage, they might still affect the accuracy during the model's training. Given that a singular value decomposition (SVD) takes place during the basis change to transform from a local to a global reduced space. The model only extracts the behavior of the training points with the dominant FOMs (i.e., design with the largest system matrix entries) and neglects the other training points. This issue is particularly true when the FOM possesses a large number of degrees of freedom (DoFs).

Therefore, we extend the idea of pBI into the basis change step to arrive at the main method of this work. This technique is the \textit{Parametric Box Reduction} (pBR) (Sect. \ref{ss:BR}) and only computes the reduced model with the training points that create a box surrounding the new evaluation point. This procedure can be executed either in the online or offline stage. It relies on box clustering to select the training points for the basis change and interpolation points, and TPWF to compute the new reduced system. Fig. \ref{fig:pBR2} depicts in \textbf{blue} the steps of this new methodology. pBR tries to predict the primary variable field with higher accuracy and a slight increase in computational effort compared to the classical approach.

\subsection{Parametric Model Order Reduction by Matrix Interpolation}

For parameter-dependent problems, Krylov subspace methods can be combined with matrix interpolation techniques to create pROM \cite{Panzer2010}. By applying the Laplace transform to the LTI system in Eq. \eqref{eq:FOM_1st} and assuming that $u(0)=0$, the transfer functions $\mathcal{H}_i$ for a first-order system is

\begin{equation}\label{eq:H1}
    \mathcal{H}_i(s)= \mathbf{D}_i\left(s\mathbf{E}_i-\mathbf{A}_i\right)^{-1} f_i\, .
\end{equation}

and repeating the procedure for the system in Eq. \eqref{eq:FOM_2nd}, we obtain the transfer function of a second-order system

\begin{equation}\label{eq:H2}
    \mathcal{H}_i(s)= \mathbf{D}_i\left(s^2\mathbf{M}_i-\mathbf{K}_i\right)^{-1} f_i\, . 
\end{equation}

The transfer function in Eq. \eqref{eq:H1} or \eqref{eq:H2} can be approximated with the help of a Krylov subspace $Kr_{r,i}$, where $r<<n$ and its basis $V_i \in  \mathbb{R}^{n \times r}$ are found with the Arnoldi process. $r$ is the number of modes or dimensions of the ROM. Its value depends on the desired accuracy. The classical algorithm of the Arnoldi method can be found in \cite{Liesen2012}; however, we implement a modified algorithm, which uses LU sparse factorization to reduce the computational effort as shown in \cite{Bechtold2005}. $V_i$ is the projection matrix to the reduced space from the full-order space. The low-order representation of the system is calculated as follows

\begin{equation}\label{eq:ReduceMatrix}
    \begin{aligned}
        \mathbf{E}_{r,i} &= {V_i}^T \mathbf{E}_i V_i \in \mathbb{R}^{r \times r}, \quad \mathbf{A}_{r,i} = {V_i}^T \mathbf{A}_i V_i \in \mathbb{R}^{r \times r}\\
        f_{r,i}  &= {V_i}^T f_i  \in \mathbb{R}^{r}, \quad \quad \quad \mathbf{D}_{r,i}  = \mathbf{D}_i V_i \in \mathbb{R}^{m \times r} \,.        
    \end{aligned}
\end{equation}

The reduced order model of the LTI can be read as: 

\begin{equation}\label{eq:ReduceSystem1}
    \begin{aligned}
        \mathbf{E}_{r,i}\dot{z}_r(t)= \mathbf{A}_{r,i}z_r(t) + f_{r,i} \,, \\
    \end{aligned}
\end{equation}

where $z$ is the reduced system variable. The same holds for second-order systems
\begin{equation}\label{eq:ReduceSystem2}
    \begin{aligned}
        \mathbf{M}_{r,i}\ddot{z}_r(t)= \mathbf{K}_{r,i}z_r(t) + f_{r,i} \,, \\
    \end{aligned}
\end{equation}
with $ \mathbf{M}_{r,i} = {V_i}^T \mathbf{M}_i V_i$ and $ \mathbf{K}_{r,i} = {V_i}^T \mathbf{K}_i V_i,$ respectively.

Assuming zero initial conditions, the ROMs in Eq. \eqref{eq:ReduceSystem1}-\eqref{eq:ReduceSystem2} can be solved. This work implements the implicit Euler method and the Newmark scheme to solve the first-order and second-order systems, respectively. The solution time is significantly reduced because of the low dimensionality of the matrices. Finally, the results are projected back to the full order space to obtain the variable field and the desired outputs

\begin{equation}\label{eq:Projection}
    \begin{aligned}
        u_i(t) &= V_iz_r(t) \, .
    \end{aligned}
\end{equation}

However, all $k$ computed local basis of the training points are not necessarily on the same reduced space; thus, a basis change from a local to a global reduced space is needed. A detailed explanation of the implemented procedure can be found in \cite{Panzer2010}. This work highlights only the essential steps.

For each training point, the basis change requires two transformation matrices, $\mathbf{T}_i$ and $\mathbf{W}_i$. These are calculated with the help of a matrix $\mathbf{R}$, which is the first $r$ columns of the left-singular vector matrix $\textbf{U}$ of the SVD of $\textbf{V}_{all}$, see Eq. \eqref{eq:RMatrix}. The $\textbf{V}_{all}$ matrix is the union of all projection matrices $V_i$. The idea behind selecting $\mathbf{R}$ in such a way is that the SVD identifies the most dominant overall behavior of the ROMs.

\begin{equation}\label{eq:RMatrix}
    \begin{aligned}
        V_{all} &= \left[{V}_{1}, \; {V}_{2}, \; \dots, {V}_{k} \right]\\
        V_{all} &= \mathbf{U} \Sigma \mathbf{V} \in \mathbb{R}^{n \times \left(r\cdot k\right)} \\
        \mathbf{R} &= U[:, 1:r] \,.\\
    \end{aligned}
\end{equation}

The transformation matrices $\textbf{W}_i$ and $\textbf{T}_i$ are calculated as follows

\begin{equation}\label{eq:TransformationMatrix}
    \begin{aligned}
        \mathbf{W}_i =  \left(V_{r,i}^T R \right)^{-1}, \quad
        \mathbf{W}_i \in \mathbb{R}^{r \times r}, \quad
        \mathbf{T}_i = R^T V_{r,i}, \quad
        \mathbf{T}_i \in \mathbb{R}^{r \times r} \,,
    \end{aligned}
\end{equation}

with these matrices, all ROMs are transformed into a global reduced space 

\begin{equation}\label{eq:pmorMatrices}
    \begin{aligned}
        \Tilde{\mathbf{M}}_{r,i} &= W_i \mathbf{M}_{r,i} {T_i}^{-1}, \quad
        &\Tilde{\mathbf{K}}_{r,i} &= W_i \mathbf{K}_{r,i} {T_i}^{-1}, \quad \\
        \Tilde{f}_{r,i} &= W_i f_{r,i}, \quad
        &\Tilde{\mathbf{D}}_{r,i} &= \mathbf{D}_{r,i}{T_i}^{-1}  \,,
    \end{aligned}
\end{equation}

with $\Tilde{\mathbf{M}}_{r,i} \in \mathbb{R}^{r \times r}$, $\Tilde{\mathbf{K}}_{r,i} \in \mathbb{R}^{r \times r}$, $\Tilde{f}_{r,i} \in \mathbb{R}^{r}$, and $\Tilde{\mathbf{D}}_{r,i} \in \mathbb{R}^{m \times r}$. The asymmetric projection in Eq. \ref{eq:pmorMatrices} could make the originally stable (symmetric positive-definite) system unstable. In \cite{geuss2016STABLEStabilityAlgorithm}, the authors proposed a Matrix Interpolation Algorithm which preserves the stability of the original system using semi-definite programming. The proposed method can be integrated with this algorithm.

The last step of the classical approach is the matrix interpolation. The primary variable $u_j(t)$ of a new parameter point $\mu_j$ is computed from the global ROMs. A new model for the additional point is calculated as a linear combination of the matrices in the global reduced space 

\begin{equation}\label{eq:MatrixInterpolationMatrices}
    \begin{aligned}
        \Tilde{\mathbf{X}}_{r,j} = \sum_{i=1}^{k}{\omega(\mu_j,\mu_i)\Tilde{\mathbf{X}}_{r,i}}\,,
    \end{aligned}
\end{equation}

where $\mathbf{X}_{r,j}$ is any of the system matrices and $ \omega$ is the weight interpolation function. The definition of this function can significantly influence the matrix interpolation results. Once the interpolated matrices are calculated, the system is solved on the reduced space and later projected into the full-order space to get the primary variable field. One of the most used functions is presented in \cite{Ananya2021}, where the Euclidean distance from the point $\mu_{j}$ to all training points is calculated ($\tilde{\omega}$). Afterward, all the weights are normalized $(\omega)$ to keep the partition of unity for the interpolation, i.e.,

\begin{equation}\label{eq:w_dis}
    \begin{aligned}
    \Tilde{\omega}(\mu_j,\mu_i) &= \frac{1}{\norm{\mu_j-\mu_i}}\\
    \omega(\mu_j,\mu_i) &= \frac{\Tilde{\omega}(\mu_j,\mu_i)}{\sum{\Tilde{\omega}(\mu_j,\mu_i)}}\,.
    \end{aligned}
\end{equation}

\subsection{Tensor Product Weight Function}\label{ss:TPWF}

\begin{figure}[ht]  
    \centering
    \resizebox{0.6\textwidth}{!}{
        \begin{tikzpicture}
            \tikzstyle{arrow} = [thick,->,>=stealth]
            \foreach \x in {0,5}  
                    \foreach \y in {0,3} 
                        \fill[tumblue] (\x, \y) circle (2pt);

            \node[font=\scriptsize,align=center] at (-0.3,0.0) {$\mu_1$};
            \node[font=\scriptsize,align=center] at (-0.3,3.0) {$\mu_3$};
            \node[font=\scriptsize,align=center] at (5.3,0.0) {$\mu_2$};
            \node[font=\scriptsize,align=center] at (5.3,3.0) {$\mu_4$};
            \node[font=\scriptsize,align=center] at (3.6,1.3) {$\mu_j$};

            \filldraw[tumorange] (3.5, 0.8) -- (3.7, .8) -- (3.6, 1) -- cycle;
            \foreach \x in {0,5}  
                    \foreach \y in {0,3}
                        \draw[arrow] (\x,\y) -- (3.6,0.9);

            \node[font=\scriptsize,align=center] at (2.00,0.75) {$\vec{d}_1$};
            \node[font=\scriptsize,align=center] at (2.00,2.15) {$\vec{d}_3$};
            \node[font=\scriptsize,align=center] at (4.7,0.65) {$\vec{d}_2$};
            \node[font=\scriptsize,align=center] at (4.7,2.00) {$\vec{d}_4$};

            \foreach \x in {8}  
                    \foreach \y in {0,3} 
                        \fill[tumblue] (\x, \y) circle (2pt);

            \node[font=\scriptsize,align=center] at (7.7,0.0) {$\mu_1$};
            \node[font=\scriptsize,align=center] at (7.7,3.0) {$\mu_3$};
            \node[font=\scriptsize,align=center] at (8.3,0.0) {$\mu_2$};
            \node[font=\scriptsize,align=center] at (8.3,3.0) {$\mu_4$};
            \node[font=\scriptsize,align=center] at (8.7,0.9) {$\mu_j$};

            \filldraw[tumorange] (7.9, 0.8) -- (8.1, .8) -- (8, 1) -- cycle;
            \foreach \x in {8}  
                    \foreach \y in {0,3}
                        \draw[arrow] (\x,\y) -- (8,0.9);

            \node[font=\scriptsize,align=center] at (7.7,0.45) {$\vec{d}_{1,y}$};
            \node[font=\scriptsize,align=center] at (7.7,2.20) {$\vec{d}_{3,y}$};
            \node[font=\scriptsize,align=center] at (8.3,0.45) {$\vec{d}_{2,y}$};
            \node[font=\scriptsize,align=center] at (8.3,2.20) {$\vec{d}_{4,y}$};

            \foreach \x in {0,5}  
                    \foreach \y in {-2} 
                        \fill[tumblue] (\x, \y) circle (2pt);

            \node[font=\scriptsize,align=center] at (0,-2.3) {$\mu_1$};
            \node[font=\scriptsize,align=center] at (0,-1.7) {$\mu_3$};
            \node[font=\scriptsize,align=center] at (5,-2.3) {$\mu_2$};
            \node[font=\scriptsize,align=center] at (5,-1.7) {$\mu_4$};
            \node[font=\scriptsize,align=center] at (3.6,-1.7) {$\mu_j$};

            \filldraw[tumorange] (3.5, -2.1) -- (3.7, -2.1) -- (3.6, -1.9) -- cycle;
            \foreach \x in {0,5}  
                    \foreach \y in {-2}
                        \draw[arrow] (\x,\y) -- (3.6,-2);

            \node[font=\scriptsize,align=center] at (2.00,-1.7) {$\vec{d}_{3,x}$};
            \node[font=\scriptsize,align=center] at (2.00,-2.3) {$\vec{d}_{1,x}$};
            \node[font=\scriptsize,align=center] at (4.25,-1.7) {$\vec{d}_{4,x}$};
            \node[font=\scriptsize,align=center] at (4.25,-2.3) {$\vec{d}_{2,x}$};
            
            \draw[->, line width=1mm, scale=2, color=tumgrey] (1.25,-.25) -- (1.25,-.5);
            \draw[->, line width=1mm, scale=2, color=tumgrey] (3,0.75) -- (3.25,0.75);

            \draw[->, line width=1mm, scale=2, color=tumgrey] (2.75,-1) -- (3,-1);
            \node[font=\small,align=center] at (6.5,-2) {$\mathbf{N_{x}}$};

            \draw[->, line width=1mm, scale=2, color=tumgrey] (3.5,-1) -- (3.75,-1);
        
            \draw[->, line width=1mm, scale=2, color=tumgrey] (4,-0.1) -- (4,-.35);
            \node[font=\small,align=center] at (8,-1) {$\mathbf{N_{y}}$};
            
            \draw[->, line width=1mm, scale=2, color=tumgrey] (4,-0.6) -- (4,-.85);
            \node[font=\normalsize,align=center] at (8,-2) {$\mathbf{\tilde{\omega}}$};

            \draw[->, line width=1mm, scale=2, color=tumgrey] (4.25,-1) -- (4.5,-1);
            \node[font=\normalsize,align=center] at (9.25,-2) {$\omega$}; 
        \end{tikzpicture}
        }
    \caption{Graphical depiction of the procedure to compute $\omega$ with tensor product weight function for a two-parameter model. $N_x$ and $N_y$ are the groupings of the local weight functions.}\label{fig:TPWF} 
\end{figure}
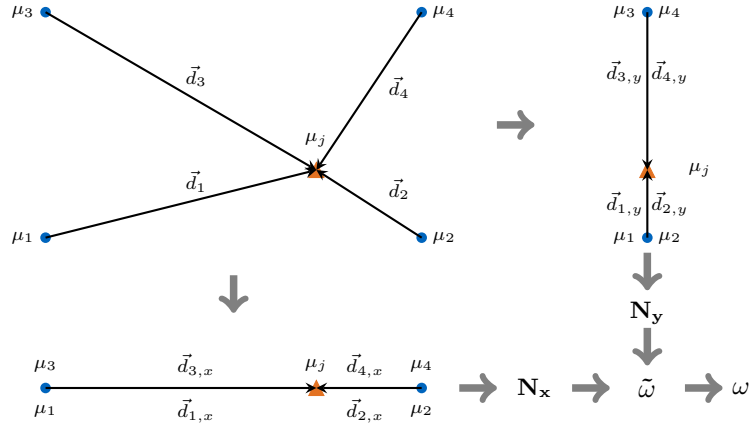

Tensor product weight function (TPWF) is a new scheme to calculate the weight function $\omega$. It is a variation of the Euclidean distance Eq. \eqref{eq:w_dis}. The main idea is to treat each parameter dimension separately and remove the need to normalize the parameter range before the interpolation step. Instead of measuring the Euclidean norm of $\vec{d}_i$ to compute $\omega$ for a new design point, the TPWF computes a local weight function $N_{i,l}$ for each dimension and training point $\mu_i \in \mathbb{T}$ 

\begin{equation}\label{eq:TPWF}
    \begin{aligned}
        \vec{d}_i &= \mu_j-\mu_i, \quad
        &\vec{d}_{i,l} &= \vec{d}_i \cdot \hat{e}_l\,,\\
        \tilde{N}_{i,l} &= \frac{1}{\norm{\vec{d}_{i,l}}}, \quad
        &N_{i,l} &= \frac{\tilde{N}_{i,l}}{\sum{\tilde{N}}_{i,l}}\,,\\
        \tilde{\omega}_{i} &= \prod^p_{l=1}{N_{i,l}}, \quad
        &\omega_i &= \frac{\tilde{\omega}_{i}}{\sum{\tilde{\omega}}_{i}}\,,\\
    \end{aligned}
\end{equation}

where $p$ is the number of design variables, $\mu_j$ is the new parameter point, and $\tilde{\omega}$ combines all the local weight functions information. To maintain the partition of unity, $\tilde{\omega}$ is normalized. Fig. \ref{fig:TPWF} exemplifies the procedure for a two-parameter model. The vectors from the training points to the new parameter point $d_i$ are projected into each dimension with a dot product. Two local weight functions $N_x$ and $N_y$ are calculated with Eq. \eqref{eq:TPWF}.

\subsection{Box Clustering}\label{ss:BC}

\begin{figure}[ht] 
    \centering  
    \begin{subfigure}[b]{0.3\textwidth} 
        \centering  
        \begin{tikzpicture}  
            \foreach \x in {0,1,2,3}
                \draw[tumgrey] (\x,0) -- (\x,3);
            \foreach \y in {0,1,2,3}  
                \draw[tumgrey] (3,\y) -- (0,\y);
            
            \foreach \x in {0,1,2,3}  
                \foreach \y in {0,1,2,3}  
                    \fill[tumblue] (\x, \y) circle (2pt);
            
            \filldraw[tumorange] (0.25, 1.25) -- (0.45, 1.25) -- (0.35, 1.45) -- cycle;
            \node[font=\scriptsize,align=center] at (0.6,1.5) {$\mu_j$}; 
        \end{tikzpicture}  
        \caption{Parameter space with training and new parameter points}\label{fig:Box1}   
    \end{subfigure}  
    \hfill  
    \begin{subfigure}[b]{0.3\textwidth}  
        \centering  
        \begin{tikzpicture}
            \foreach \x in {0,1,2,3}
                \draw[tumgrey] (\x,0) -- (\x,3);
            \foreach \y in {0,1,2,3}  
                \draw[tumgrey] (3,\y) -- (0,\y);
                
            \foreach \x in {0,1,2,3}  
                \foreach \y in {0,1,2,3}  
                    \fill[tumblue] (\x, \y) circle (2pt);
            \foreach \x in {0}  
                \foreach \y in {0,1,2,3}  
                    \node[star, fill=tumgreen, minimum size=2pt, scale = 0.35,  star point height=.25cm] at (\x, \y) {};
            \filldraw[tumorange] (0.25, 1.25) -- (0.45, 1.25) -- (0.35, 1.45) -- cycle;
            \node[font=\scriptsize,align=center] at (0.6,1.5) {$\mu_j$}; 
        \end{tikzpicture}  
        \caption{Nearest smaller points from $\mu_j$ in dimension 1}\label{fig:Box2}
    \end{subfigure}  
    \hfill  
    \begin{subfigure}[b]{0.3\textwidth}  
        \centering  
        \begin{tikzpicture}
            \foreach \x in {0,1,2,3}
                \draw[tumgrey] (\x,0) -- (\x,3);
            \foreach \y in {0,1,2,3}  
                \draw[tumgrey] (3,\y) -- (0,\y);
                
            \foreach \x in {0,1,2,3}  
                \foreach \y in {0,1,2,3}  
                    \fill[tumblue] (\x, \y) circle (2pt);
            \foreach \x in {1}  
                \foreach \y in {0,1,2,3}
                    \node[star, fill=tumgreen, minimum size=2pt, scale = 0.35,  star point height=.25cm] at (\x, \y) {};
            \filldraw[tumorange] (0.25, 1.25) -- (0.45, 1.25) -- (0.35, 1.45) -- cycle;
            \node[font=\scriptsize,align=center] at (0.6,1.5) {$\mu_j$}; 
        \end{tikzpicture}  
        \caption{Nearest larger points from $\mu_j$ in dimension 1}\label{fig:Box3}
    \end{subfigure}  
     \vskip\baselineskip
     \centering  
    \begin{subfigure}[b]{0.3\textwidth} 
        \centering  
        \begin{tikzpicture}

            \foreach \x in {0,1,2,3}
                \draw[tumgrey] (\x,0) -- (\x,3);
            \foreach \y in {0,1,2,3}  
                \draw[tumgrey] (3,\y) -- (0,\y);
            \foreach \x in {0,1,2,3}  
                \foreach \y in {0,1,2,3}  
                    \fill[tumblue] (\x, \y) circle (2pt);
            \foreach \x in {0,1,2,3}  
                \foreach \y in {1}
                    \node[star, fill=tumgreen, minimum size=2pt, scale = 0.35,  star point height=.25cm] at (\x, \y) {};
            \filldraw[tumorange] (0.25, 1.25) -- (0.45, 1.25) -- (0.35, 1.45) -- cycle;
            \node[font=\scriptsize,align=center] at (0.6,1.5) {$\mu_j$}; 
            
        \end{tikzpicture}  
        \caption{Nearest smaller points from $\mu_j$ in dimension 2}\label{fig:Box4}  
    \end{subfigure}  
    \hfill  
    \begin{subfigure}[b]{0.3\textwidth}  
        \centering  
        \begin{tikzpicture}

            \foreach \x in {0,1,2,3}
                \draw[tumgrey] (\x,0) -- (\x,3);
            \foreach \y in {0,1,2,3}  
                \draw[tumgrey] (3,\y) -- (0,\y);
                
            \foreach \x in {0,1,2,3}  
                \foreach \y in {0,1,2,3}  
                    \fill[tumblue] (\x, \y) circle (2pt);
            \foreach \x in {0,1,2,3} 
                \foreach \y in {2}
                    \node[star, fill=tumgreen, minimum size=2pt, scale = 0.35,  star point height=.25cm] at (\x, \y) {};
            \filldraw[tumorange] (0.25, 1.25) -- (0.45, 1.25) -- (0.35, 1.45) -- cycle;
            \node[font=\scriptsize,align=center] at (0.6,1.5) {$\mu_j$}; 
        \end{tikzpicture}  
        \caption{Nearest larger points from $\mu_j$ in dimension 2}\label{fig:Box5}  
    \end{subfigure}  
    \hfill  
    \begin{subfigure}[b]{0.3\textwidth}   
        \centering  
        \begin{tikzpicture}

            \foreach \x in {0,1,2,3}
                \draw[tumgrey] (\x,0) -- (\x,3);
            \foreach \y in {0,1,2,3}  
                \draw[tumgrey] (3,\y) -- (0,\y);
                
            \foreach \x in {0,1,2,3}  
                \foreach \y in {0,1,2,3}  
                    \fill[tumblue] (\x, \y) circle (2pt);
            \foreach \x in {0,1}  
                \foreach \y in {1,2}
                    \node[star, fill=tumgreen, minimum size=2pt, scale = 0.35,  star point height=.25cm] at (\x, \y) {};
            \filldraw[tumorange] (0.25, 1.25) -- (0.45, 1.25) -- (0.35, 1.45) -- cycle;
            \node[font=\scriptsize,align=center] at (0.6,1.5) {$\mu_j$}; 
        \end{tikzpicture}  
        \caption{Training points selected by box interpolation for $\mu_j$}\label{fig:Box6}    
    \end{subfigure}
    \caption{Workflow schematic of the parametric point selection methodology with box clustering. Blue circles are the training points, and green stars are the selected points.}\label{fig:Box} 
\end{figure}
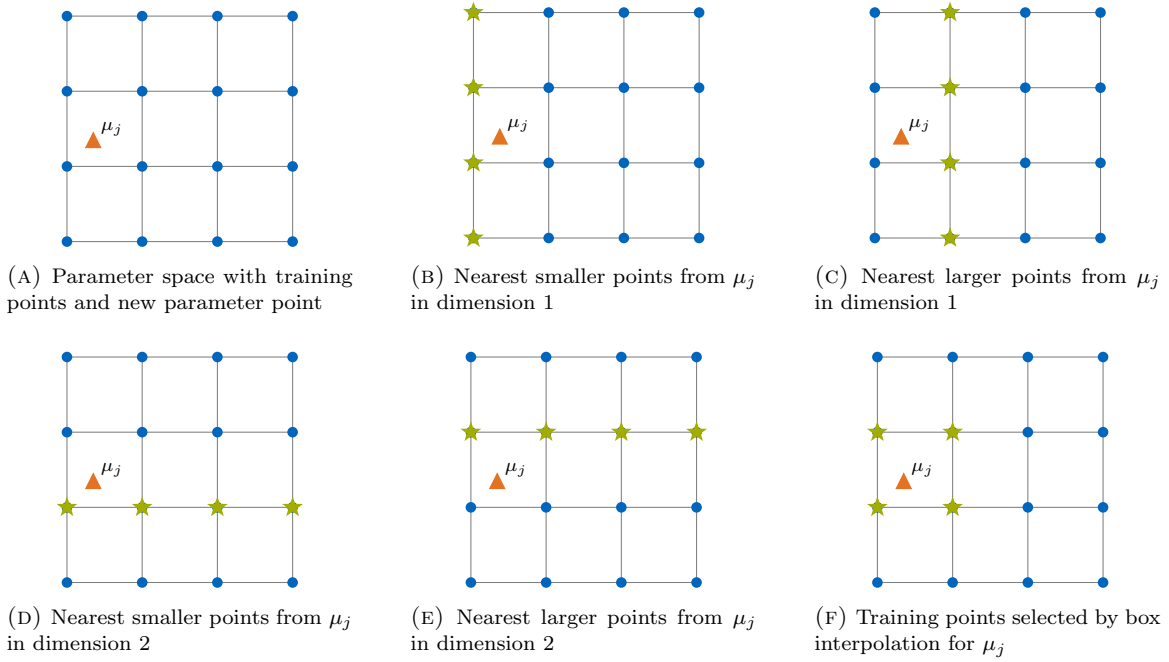

Box clustering takes inspiration from the nearest neighbor algorithm. However, it does not depend on a prior normalization of the parameter space. It finds the training points that create a box surrounding $\mu_j$ instead of the nearest ones. For this first formulation, the training points must be in a structured grid over the parameter space. In the case of a few design parameters, the trained model does not suffer from the dimensionality curse. However, the required training points grow exponentially as the parameter space dimension increases. In particular, for a parameter space with $p$ dimensions, the algorithm finds the $2^p$ training points that create the box around the new parameter point $\mu_j$. The sampling process is relatively simple since the points are in a structured grid. Figure \ref{fig:Box} depicts the main ideas of the algorithm for a two-parameter design space; similar to TPWF, the operations are executed for each parameter and then combined. In each dimension, the nearest smaller and nearest larger values are identified. The nearest smaller values of the first parameter are illustrated in Fig. \ref{fig:Box2}. Similarly, Figure \ref{fig:Box3} highlights the nearest highest for the first parameter. The selected training points must fulfill the following condition: all their design values must be either a nearest lowest or a nearest highest point from $\mu_j$ for all dimensions. The selected points are displayed in Figure \ref{fig:Box6}. Note that a formulation for higher dimensions with unstructured grids for box clustering is still an open question; ideas from \cite{Geuss2014, Geuss2016,Brunsch2017} for matrix interpolation in sparse grids can be employed. 

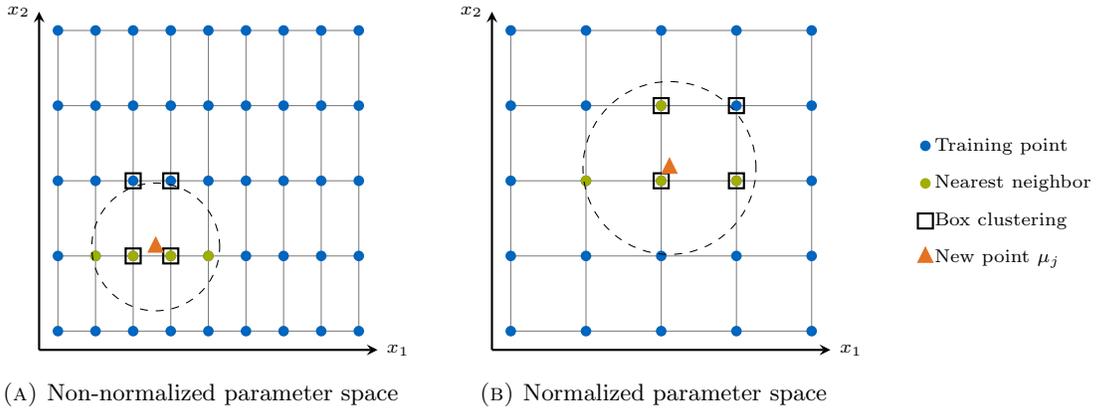
\begin{figure}[ht]  
    \centering   
    \begin{subfigure}[b]{0.35\textwidth}  
        \centering  
        \begin{tikzpicture}
            \tikzstyle{arrow} = [thick,->,>=stealth]
            \draw[arrow] (-0.25,-0.25) -- (-0.25,4.25);
            \draw[arrow] (-0.25,-0.25) -- (4.25,-0.25);
            \node[left, font=\scriptsize] at (-0.25,4.25) {$x_2$};
            \node[right, font=\scriptsize] at (4.25,-0.25) {$x_1$};
            
            \foreach \x in {0,0.5,1,1.5,2,2.5,3,3.5,4} 
                \draw[tumgrey] (\x,0) -- (\x,4);
            \foreach \y in {0,1,2,3,4}   
                \draw[tumgrey] (0,\y) -- (4,\y);
            
            \foreach \x in {0,0.5,1,1.5,2,2.5,3,3.5,4}  
                \foreach \y in {0,1,2,3,4}  
                    \fill[tumblue] (\x, \y) circle (2pt);
            \foreach \x in {0.5,1,1.5,2}  
                \foreach \y in {1}  
                    \fill[tumgreen] (\x, \y) circle (2pt);
            \foreach \x in {1,1.5}  
                \foreach \y in {1,2}  
                    \draw[black, thick] (\x-0.1, \y-0.1) rectangle (\x+0.1, \y+0.1);
            \filldraw[tumorange] (1.2, 1.05) -- (1.4, 1.05) -- (1.3, 1.25) -- cycle;
            \draw[dashed] (1.3,1.12) circle (0.85);
        \end{tikzpicture}  
        \caption{Non-normalized parameter space}\label{fig:NearNorm}  
    \end{subfigure}  
    \hfill   
    \begin{subfigure}[b]{0.35\textwidth}  
        \centering  
        \begin{tikzpicture}  
            \tikzstyle{arrow} = [thick,->,>=stealth]
            \draw[arrow] (-0.25,-0.25) -- (-0.25,4.25);
            \draw[arrow] (-0.25,-0.25) -- (4.25,-0.25);
            \node[left, font=\scriptsize] at (-0.25,4.25) {$x_2$};
            \node[right, font=\scriptsize] at (4.25,-0.25) {$x_1$};
            \foreach \x in {0,1,2,3,4} 
                \draw[tumgrey] (\x,0) -- (\x,4);
            \foreach \y in {0,1,2,3,4}   
                \draw[tumgrey] (0,\y) -- (4,\y);
                
            \foreach \x in {0,1,2,3,4}  
                \foreach \y in {0,1,2,3,4}  
                    \fill[tumblue] (\x, \y) circle (2pt);
            \foreach \x in {1,2,3}  
                \foreach \y in {2}  
                    \fill[tumgreen] (\x, \y) circle (2pt);
            \fill[tumgreen] (2, 3) circle (2pt);
            \foreach \x in {2,3}  
                \foreach \y in {2,3}  
                    \draw[black, thick] (\x-0.1, \y-0.1) rectangle (\x+0.1, \y+0.1);
            \filldraw[tumorange] (2.01, 2.1) -- (2.21, 2.1) -- (2.11, 2.3) -- cycle;
            \draw[dashed] (2.11,2.17) circle (1.15);
        \end{tikzpicture}
        \caption{Normalized parameter space}\label{fig:NearPoints}
    \end{subfigure}
    \begin{subfigure}[b]{0.25\textwidth}  
        \centering  
        \begin{tikzpicture}   
            \fill[tumblue] (5, 4.5) circle (2pt);
            \node[right,font=\fontsize{7}{10}\selectfont] at (5,4.5) {Training point};
            \fill[tumgreen] (5, 4) circle (2pt);
            \node[right,font=\fontsize{7}{10}\selectfont] at (5,4) {Nearest neighbor};
            \draw[black, thick] (5-0.1, 3.5-0.1) rectangle (5+0.1, 3.5+0.1);
            \node[right,font=\fontsize{7}{10}\selectfont] at (5,3.5) {Box clustering};
            \filldraw[tumorange] (4.9, 2.95) -- (5.1, 2.95) -- (5, 3.15) -- cycle;
            \node[right,font=\fontsize{7}{10}\selectfont] at (5,3) {New point $\mu_j$};
            \fill[tumblue] (5, 1) circle (0.01pt);
        \end{tikzpicture} 
    \end{subfigure}
    \caption{Example illustration depicting nearest neighbor challenges in normalized and non-normalized parameter space. Identification of the kNN with the dashed circle. Box clustering points are not affected by the normalization of the parameter space.}\label{fig:NearBI} 
\end{figure}

Although nearest neighbor is widely employed, it presents some issues as depicted in Fig. \ref{fig:NearBI}. Fig. \ref{fig:NearNorm} illustrates how important normalization is. It shows a non-normalized parameter grid, where the radius of the nearest neighbor only selects points with the same value for $x_2$. Therefore, the new reduced model behaves as if the new point $\mu_j$ has the same $x_2$ as the nearest neighbor training points, even though the real $x_2$ value is higher, thus artificially influencing the pROM behavior. In Fig. \ref{fig:NearPoints}, a MinMax normalization of the parameter range is depicted; however, the method leads to misleading results for a new simulation near one of the training points since the nearest neighbor introduces an artificial influence on the interpolation. Three of the four points have the same $x_2$ value. Therefore, the design parameter corresponding to $x_2$  has a more significant impact on the trained model than the parameter $x_1$. By applying box clustering, a consistent and meaningful selection of training points (i.e., box points) is achieved without the normalization of the parameter space. Therefore, solving the issues from k's nearest neighbor highlighted in Fig. \ref{fig:NearBI}.

\subsection{Parametric Box Interpolation}\label{ss:BI}

In \cite{Yuan2021b}, the authors chose the kNN strategy to select the ROMs for the matrix interpolation. However, as illustrated in Fig. \ref{fig:NearBI}, kNN has drawbacks and depends on the normalization step. Hence, we introduce the pBI method, which computes the matrix interpolation using the TPWF with the points selected with the box clustering scheme. 

\begin{figure}[ht]
    \centering
        \begin{subfigure}[b]{0.223\textwidth}
            \centering
            \resizebox{0.95\textwidth}{!}{
            \begin{tikzpicture}
                \foreach \x in {0,2} 
                    \draw[tumgrey] (\x,0) -- (\x,2);
                \foreach \y in {0,2}   
                    \draw[tumgrey] (0,\y) -- (2,\y);
                
                \foreach \x in {0,2}  
                    \foreach \y in {0,2}  
                        \fill[tumblue] (\x, \y) circle (2pt);
                \filldraw[tumorange] (0.8, 0.2) -- (0.6, 0.2) -- (0.7, 0.4) -- cycle;
                \node[scale=0.85] at (0.5,0.4) {\tiny $\mu_j$};
                
                \foreach \x in {1,3}
                    \node[scale=0.85] at (\x-1,-0.2) {\tiny$\mu_\x$};
                
                \foreach \x in {7,9}
                    \node[scale=0.85] at (\x-7,2.2) {\tiny $\mu_\x$};
            \end{tikzpicture}
            }
            \caption{Four training points and low value of design variables}
            \label{fig:pBI_k4L}
        \end{subfigure}
        \hfill
        \begin{subfigure}[b]{0.223\textwidth}
            \centering
            \resizebox{0.95\textwidth}{!}{
            \begin{tikzpicture}
                \foreach \x in {0,1,2} 
                    \draw[tumgrey] (\x,0) -- (\x,2);
                \foreach \y in {0,1,2}   
                    \draw[tumgrey] (0,\y) -- (2,\y);
                
                \foreach \x in {0,1,2}  
                    \foreach \y in {0,1,2}  
                        \fill[tumblue] (\x, \y) circle (2pt);
                \filldraw[tumorange] (0.8, 0.2) -- (0.6, 0.2) -- (0.7, 0.4) -- cycle;
                \node[scale=0.85] at (0.5,0.4) {\tiny $\mu_j$};
                
                \tikzstyle{arrow} = [->,>=stealth]
                \draw[arrow] (2,2) -- (0.7,0.28);
                \draw[arrow] (0,0) -- (0.7,0.28);
                
                \foreach \x in {1,2,3}
                    \node[scale=0.85] at (\x-1,-0.2) {\tiny$\mu_\x$};
                \foreach \x in {5,6}
                    \node[scale=0.85] at (\x-4.15,1.2) {\tiny $\mu_\x$};
                \foreach \x in {4}
                    \node[scale=0.85] at (\x-3.8,1.2) {\tiny $\mu_\x$};

                \foreach \x in {7,8,9}
                    \node[scale=0.85] at (\x-7,2.2) {\tiny $\mu_\x$};
            \end{tikzpicture}
            }
            \caption{Nine training points and low value of design variables}
            \label{fig:pBI_k9L}
        \end{subfigure}
        \hfill
        \begin{subfigure}[b]{0.223\textwidth}
            \centering
            \resizebox{0.95\textwidth}{!}{
            \begin{tikzpicture}
                \foreach \x in {0,2} 
                    \draw[tumgrey] (\x,0) -- (\x,2);
                \foreach \y in {0,2}   
                    \draw[tumgrey] (0,\y) -- (2,\y);
                
                \foreach \x in {0,2}  
                    \foreach \y in {0,2}  
                        \fill[tumblue] (\x, \y) circle (2pt);
                \filldraw[tumorange] (1.8, 1.6) -- (1.6, 1.6) -- (1.7, 1.8) -- cycle;
                \node[scale=0.85] at (1.5,1.75) {\tiny $\mu_j$};
                
                \foreach \x in {1,3}
                    \node[scale=0.85] at (\x-1,-0.2) {\tiny$\mu_\x$};
                
                \foreach \x in {7,9}
                    \node[scale=0.85] at (\x-7,2.2) {\tiny $\mu_\x$};
            \end{tikzpicture}
            }
            \caption{Four training points, high value of design variables}
            \label{fig:pBI_k4H}
        \end{subfigure}
        \hfill
        \begin{subfigure}[b]{0.223\textwidth}
            \centering
            \resizebox{0.95\textwidth}{!}{
            \begin{tikzpicture}
                \foreach \x in {0,1,2} 
                    \draw[tumgrey] (\x,0) -- (\x,2);
                \foreach \y in {0,1,2}   
                    \draw[tumgrey] (0,\y) -- (2,\y);
                
                \foreach \x in {0,1,2}  
                    \foreach \y in {0,1,2}  
                        \fill[tumblue] (\x, \y) circle (2pt);
                \filldraw[tumorange] (1.8, 1.6) -- (1.6, 1.6) -- (1.7, 1.8) -- cycle;
                \node[scale=0.85] at (1.5,1.75) {\tiny $\mu_j$};
                
                \tikzstyle{arrow} = [->,>=stealth]
                \draw[arrow] (2,2) -- (1.7,1.68);
                \draw[arrow] (0,0) -- (1.7,1.68);
                
                \foreach \x in {1,2,3}
                    \node[scale=0.85] at (\x-1,-0.2) {\tiny$\mu_\x$};
                \foreach \x in {5,6}
                    \node[scale=0.85] at (\x-4.15,1.2) {\tiny $\mu_\x$};
                \foreach \x in {4}
                    \node[scale=0.85] at (\x-3.8,1.2) {\tiny $\mu_\x$};

                \foreach \x in {7,8,9}
                    \node[scale=0.85] at (\x-7,2.2) {\tiny $\mu_\x$};
            \end{tikzpicture}
            }
            \caption{Nine training points, high value of design variables}
            \label{fig:pBI_k9H}
        \end{subfigure}
    \caption{Parameter space for two parameter points with low and high-value design parameters. Four and nine training points are illustrated for each case.}
    \label{fig:pBI}
\end{figure}
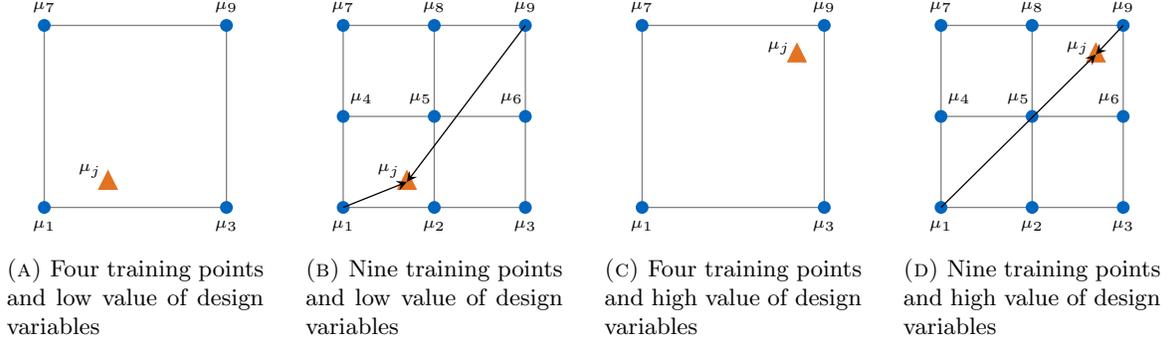

The motivation to not compute the new pROM with all the training points is illustrated in Fig. \ref{fig:pBI}. In Fig. \ref{fig:pBI_k4L}, only four training points are selected in the parameter space to predict the behavior of a new parameter point $\mu_j$ with low values of the design variables. To increase the accuracy, more points are selected for training the model, Fig. \ref{fig:pBI_k9L}. The training point $\mu_1$ is expected to have a greater influence than $\mu_9$ on the reduced model because it holds $\omega_1 > \omega_9$. Given the large parameter ranges, the behavior of the system differs largely between training points. This is reflected in the magnitude of the system matrices, e.g., an increase in Young's modulus increases the stiffness matrix $\mathbf{K}$ entries. Thus the system matrix entries of training point $\mu_9$ are much larger than the entries from $\mu_1$, which is reflected in the reduced matrices $\mathbf{X}_{r,i}$. Therefore $ \omega_1\textbf{X}_{r,1} \ll \omega_9 \textbf{X}_{r,9}$, given $\mu_9$ a higher impact on the reduced model. This influence comes from the matrix interpolation procedure and not from the physics of the problem. By using only the box points, this artificial influence is neglected. Similarly, Fig. \ref{fig:pBI_k4H} and Fig. \ref{fig:pBI_k9H} depict a new point with high design parameter values. Compared to the previous case, pBI enhances or, at a minimum, preserves the performance when increasing the number of training points because the new parameter point $\mu_j$ is near the dominant training point $\mu_9$, i.e., the artificial influence from $\mu_1$, a distant training point, is negligible. 

\subsection{Parametric Box Reduction}\label{ss:BR}

The parametric Box Reduction (pBR) extends the idea from pBI into the basis change stage. In \cite{Panzer2010}, the authors introduce a similar concept, where for each new parameter evaluation, they multiply $V_i$ with the weight interpolation function to construct a modified $\bar{\mathbf{V}}_{all} \in \mathbb{R}^{n \times \left(r\cdot k\right)}$, which is used for the SVD. In this work, instead of combining all $k$ projection matrices with the respective weight into $\bar{\mathbf{V}}_{all}$, only the $2^p$ box points create the modified $\tilde{\mathbf{V}}_{all} \in \mathbb{R}^{n \times \left(r\cdot 2^p\right)}$. The computational cost of the weighted method in \cite{Panzer2010} is one of the main drawbacks, especially when the number of training points $k$ increases. The SVD complexity \cite{li2019tutorialcomplexityanalysissingular} for a matrix $A \in \mathbb{R}^{n\times m}$ with $m\ll n$ is $\mathcal{O}(4n^2m^2)$. Thus, the pBR has a speed-up of $(k/2^p)^2$ with respect to the weighted matrix interpolation in \cite{Panzer2010}. Compared to the classical approach \cite{Panzer2010}, the computational time for pBR increases because of the need to compute an SVD at each new parameter point. Nonetheless, this computational time can be reduced and, in some cases, perform better than the nonweighted matrix interpolation if the SVD results are stored for all the possible boxes of the parameter space in a pre-processing step. However, in that case, the memory requirement increases.

During the basis change, the SVD identifies the most important information from the  $\textbf{V}_{all}$ matrix to compute the transformation matrices ($\textbf{T}_i$ and $\textbf{W}_i$). These matrices change the ROMs from a local to a global reduced space that maintains the main characteristic of the HFM. Since the SVD extracts the most energetic modes, relevant information on some training points can be lost while simultaneously amplifying the influence of less pertinent points. This occurs mainly when the magnitude of the HFM matrices entries significantly differ between training points. The parameter points with larger entries have more influence on the basis change, while the small ones are neglected, spuriously influencing the pMOR results. pBR aims to reduce the non-physical influence of farther-away training points, especially where wide parameter ranges are used. The weighted approach proposed by \cite{Panzer2010} can also reduce this influence. However, the issue depicted in \ref{fig:pBI}, where $ \omega_1V_{1} \ll \omega_9 V_{9}$ could still occur when more training points are used.

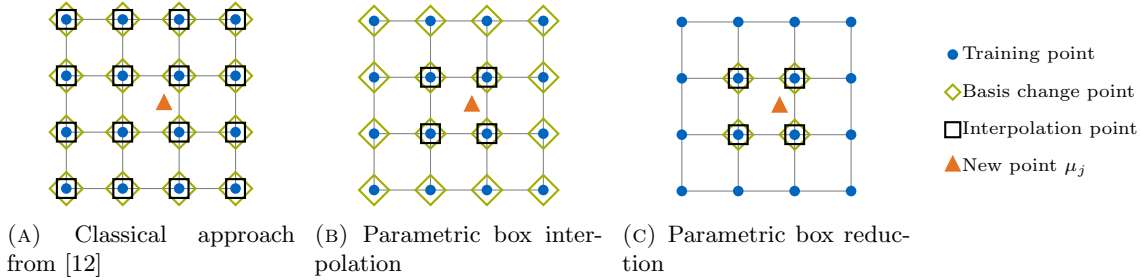
\begin{figure}[ht]  
    \centering   
    \begin{subfigure}[b]{0.23\textwidth}  
        \centering
        \resizebox{0.95\textwidth}{!}{
        \begin{tikzpicture}

            \foreach \x in {0,0.75,1.5,2.25} 
                \draw[tumgrey] (\x,0) -- (\x,2.25);
            \foreach \y in {0,0.75,1.5,2.25}   
                \draw[tumgrey] (0,\y) -- (2.25,\y);
                
            \foreach \x in {0,0.75,1.5,2.25}  
                \foreach \y in {0,0.75,1.5,2.25}  
                    \fill[tumblue] (\x, \y) circle (2pt);
            \foreach \x in {0,0.75,1.5,2.25}  
                \foreach \y in {0,0.75,1.5,2.25}  
                    \draw[tumgreen, thick] (\x+0.22, \y) -- (\x, \y-0.22) -- (\x-0.22, \y) -- (\x, \y+0.22) --cycle; 
            \foreach \x in {0,0.75,1.5,2.25}  
                \foreach \y in {0,0.75,1.5,2.25}  
                    \draw[black, thick] (\x-0.125, \y-0.125) rectangle (\x+0.125, \y+0.125);
            \filldraw[tumorange] (1.2, 1.05) -- (1.4, 1.05) -- (1.3, 1.25) -- cycle;
        \end{tikzpicture}
        }
        \caption{Classical approach from \cite{Panzer2010}}\label{fig:Classical}  
    \end{subfigure}  
    \hfill   
    \begin{subfigure}[b]{0.23\textwidth}  
        \centering
        \resizebox{0.95\textwidth}{!}{
        \begin{tikzpicture}

            \foreach \x in {0,0.75,1.5,2.25} 
                \draw[tumgrey] (\x,0) -- (\x,2.25);
            \foreach \y in {0,0.75,1.5,2.25}   
                \draw[tumgrey] (0,\y) -- (2.25,\y);
                
            \foreach \x in {0,0.75,1.5,2.25}  
                \foreach \y in {0,0.75,1.5,2.25}  
                    \fill[tumblue] (\x, \y) circle (2pt);
            \foreach \x in {0,0.75,1.5,2.25}  
                \foreach \y in {0,0.75,1.5,2.25}  
                    \draw[tumgreen, thick] (\x+0.2, \y) -- (\x, \y-0.2) -- (\x-0.2, \y) -- (\x, \y+0.2) --cycle; 
            \foreach \x in {0.75,1.5}  
                \foreach \y in {0.75,1.5}  
                    \draw[black, thick] (\x-0.125, \y-0.125) rectangle (\x+0.125, \y+0.125);
            \filldraw[tumorange] (1.2, 1.05) -- (1.4, 1.05) -- (1.3, 1.25) -- cycle;
        \end{tikzpicture}
        }
        \caption{Parametric box interpolation}\label{fig:BIPoints}
    \end{subfigure}
    \hfill   
    \begin{subfigure}[b]{0.23\textwidth}  
        \centering
        \resizebox{0.86\textwidth}{!}{
        \begin{tikzpicture}
            \foreach \x in {0,0.75,1.5,2.25} 
                \draw[tumgrey] (\x,0) -- (\x,2.25);
            \foreach \y in {0,0.75,1.5,2.25}   
                \draw[tumgrey] (0,\y) -- (2.25,\y);
                
            \foreach \x in {0,0.75,1.5,2.25}  
                \foreach \y in {0,0.75,1.5,2.25}  
                    \fill[tumblue] (\x, \y) circle (2pt);
            \foreach \x in {0.75,1.5}  
                \foreach \y in {0.75,1.5}  
                    \draw[tumgreen, thick] (\x+0.2, \y) -- (\x, \y-0.2) -- (\x-0.2, \y) -- (\x, \y+0.2) --cycle;
            \foreach \x in {0.75,1.5}  
                \foreach \y in {0.75,1.5}  
                    \draw[black, thick] (\x-0.125, \y-0.125) rectangle (\x+0.125, \y+0.125);
            \filldraw[tumorange] (1.2, 1.05) -- (1.4, 1.05) -- (1.3, 1.25) -- cycle;
            \fill[tumblue] (0, -0.2) circle (0.01pt);
            \fill[tumblue] (0, 3.2) circle (0.01pt);
        \end{tikzpicture}
        }
        \caption{Parametric box reduction}\label{fig:BRPoints}
    \end{subfigure}
    \hfill 
    \begin{subfigure}[b]{0.2\textwidth}  
        \centering  
        \begin{tikzpicture}   
            \fill[tumblue] (5, 4.5) circle (2pt);
            \node[right,font=\fontsize{7}{10}\selectfont] at (5,4.5) {Training point};
            \draw[tumgreen, thick] (5+0.1, 4) -- (5, 
            4-0.1) -- (5-0.1, 4) -- (5, 4+0.1) --cycle;
            \node[right,font=\fontsize{7}{10}\selectfont] at (5,4) {Basis change point};
            \draw[black, thick] (5-0.1, 3.5-0.1) rectangle (5+0.1, 3.5+0.1);
            \node[right,font=\fontsize{7}{10}\selectfont] at (5,3.5) {Interpolation point};
            \filldraw[tumorange] (4.9, 2.95) -- (5.1, 2.95) -- (5, 3.15) -- cycle;
            \node[right,font=\fontsize{7}{10}\selectfont] at (5,3) {New point $\mu_j$};
            \fill[tumblue] (5, 1.5) circle (0.01pt);
        \end{tikzpicture}
    \end{subfigure}
    \caption{Visualization of the selected points for basis change and matrix interpolation in the different model order reduction strategies.}\label{fig:BoxRed} 
\end{figure}

Figure \ref{fig:BoxRed} highlights three reduction schemes to better illustrate the different methods. The classical approach \cite{Panzer2010} is depicted in Fig. \ref{fig:Classical}, where all the training points are used in the basis change and matrix interpolation stage. In contrast, Fig. \ref{fig:BIPoints} illustrates the pBI method; only the points surrounding $\mu_j$ are sampled for the interpolation, but all the training points construct the $V_{all}$ matrix. Figure \ref{fig:BRPoints} shows that pBR computes the global space with only the box points.

\section{Results \& Discussion}\label{ss:results}

In this section, we showcase the developed techniques by two examples from different domains. The first test case is a three-dimensional cantilever Timoshenko beam introduced in \cite{Panzer2009}, commonly used as an academic benchmark \cite{Schopper2024}. In contrast, the second example is an industrial application case, a transient thermal simulation of a power module. For both cases, a two-parameter model is generated with the classical approach, parametric box interpolation, and parametric box reduction.

In the evaluation, a prime focus is set on the effect of varying numbers of training simulations $k$ in the training phase. A higher number of training points often leads to improved model accuracy but compromises the overall cost efficiency. Hence, we investigate the influence of $k$ on the performance of the pROM. For both application cases, the design parameters of the training set are defined on a structure grid (see Section \ref{ss:BC}). Additionally, the validation points are pseudorandomly selected with a Latin hypercube sampling strategy to explore the entire parameter space with a small number of sample points. 

The performance in terms of accuracy and computational cost is measured with various indicators. These are calculated for each reduced model to compare the efficiency of the presented methods. To evaluate the method's accuracy across the entire parametric space, full error surface responses for the cantilever beam are employed, and bubble charts for the power module, where different disks are plotted. Each disk represents one model evaluation with its parameter configuration according to the x- and y-coordinate. The size of the disk indicates the magnitude of the error. In addition, multiple Figures show the evolution of the error in selected validation points for a quantitative comparison. For performance assessment, CPU time plots compare the total computation time for a FOM and a pROM. Furthermore, an analysis identifies the \textit{break-even point}, the number of evaluations with an equivalent computational time of FOM and pROM. For evaluations conducted beyond this break-even point, the pROM demonstrates superior computational efficiency, including the construction phase. 

\subsection{Cantilever Timoshenko beam}

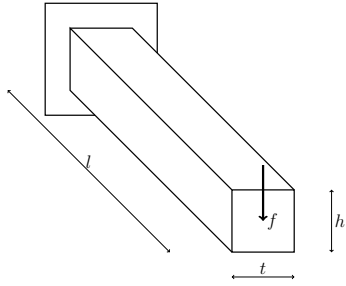
\begin{figure}[ht]
    \centering
    \resizebox{0.25\textwidth}{!}{
    \begin{tikzpicture}
        \draw[line width=0.5mm] (7.5,0) rectangle (10,2.5);
        \draw[line width=0.5mm] (7.5,0)-- (1,6.5) -- (1,9) -- (7.5,2.5);
        \draw[line width=0.5mm] (10,2.5)-- (3.5,9) -- (1,9);
        \draw[line width=0.5mm] (2,5.5) -- (0,5.5) -- (0,10) -- (4.5, 10) -- (4.5,8);
        \draw[line width=0.25mm, <->] (11.5,0) -- (11.5, 2.5) node[midway, right] {\Huge$h$};
        \draw[line width=0.25mm, <->] (7.5,-1) -- (10, -1) node[midway, above] {\Huge$t$};
         \draw[line width=0.25mm, <->] (-1.5,6.5) -- (5, 0) node[midway, above] {\Huge$l$};
         \draw[line width=0.9mm, ->] (8.75,3.5) -- (8.75, 1.25) node[above, right] {\Huge$f$};
    \end{tikzpicture}
    }
    \caption{Three-dimensional cantilever Timoshenko beam.}
    \label{fig:Beam}
\end{figure}

Figure \ref{fig:Beam} depicts the 3D beam model \cite{Panzer2009}. It is clamped on one end, and a transient load $f$ is applied to the other end. We consider an undamped system for this work.

The characteristic equation can be written as a second-order LTI system, such as Eq. \eqref{eq:FOM_2nd}. The beam is discretized into 100 beam elements, every node has six degrees of freedom: three translational displacements and three rotational degrees of freedom  (i.e., $n=600$). Its time step size is $dt = 10^{-6}$ and the solving time $t_{FOM} = 5.72\times 10^{-4}\;\mathrm{s}$ for each iteration of the Newmark scheme. We choose the Normalized Root Mean Square Error as an accuracy indicator.

\begin{equation}\label{eq:NRMSE}
    NRMSE = \frac{\sqrt{\frac{\sum\limits_{i=1}^{n}{\left(u_s- u_{r,i}\right)^2}}{n}}}{u_{s,tip}} \,,
\end{equation}

where $u_s$ is the HFS solution field, $u_{r,i}$ the predicted field, and $u_{tip}$ the displacement on the tip. The mass density $\rho$ and Young's modulus $E$ are chosen as design parameters for the numerical experiments. The specific parameters describing the simulation set-up are indicated in Table \ref{tab:B_2p}. 

\begin{table}[!ht]
    \centering
    \caption{Simulation set-up of three-dimensional cantilever Timoshenko beam with geometry and material properties.}
    {\footnotesize
    \begin{tabular}{c|c|c}
         \rule{0pt}{1em}Parameter&Range/Value&Unit\\  \hline \rule{0pt}{1em}
            Length $l$ & $1.0$ & $\mathrm{m}$ \\ \rule{0pt}{1em}
            Thickness $t$ & $0.01$ & $\mathrm{m}$ \\ \rule{0pt}{1em}
            Height $h$ & $0.01$ & $\mathrm{m}$ \\ \rule{0pt}{1em}
            Mass density $\rho$ & $[6350, 8850]$ & $\mathrm{kg}/\mathrm{m}^3$\\ \rule{0pt}{1em}
            Young's modulus $E$ & $[10,410]$ & $\mathrm{GPa}$\\ \rule{0pt}{1em}
            Poisson's ratio $\nu$ & $0.3$ & -\\ \rule{0pt}{1em}
            Simulation time $t_{end}$ & $0.01$ & $\mathrm{s}$ \\
    \end{tabular}
    }
    \label{tab:B_2p}
\end{table}

\begin{figure}[ht]
    \centering
    \begin{subfigure}[b]{0.3\textwidth}
        \centering
        \includegraphics[width=\textwidth]{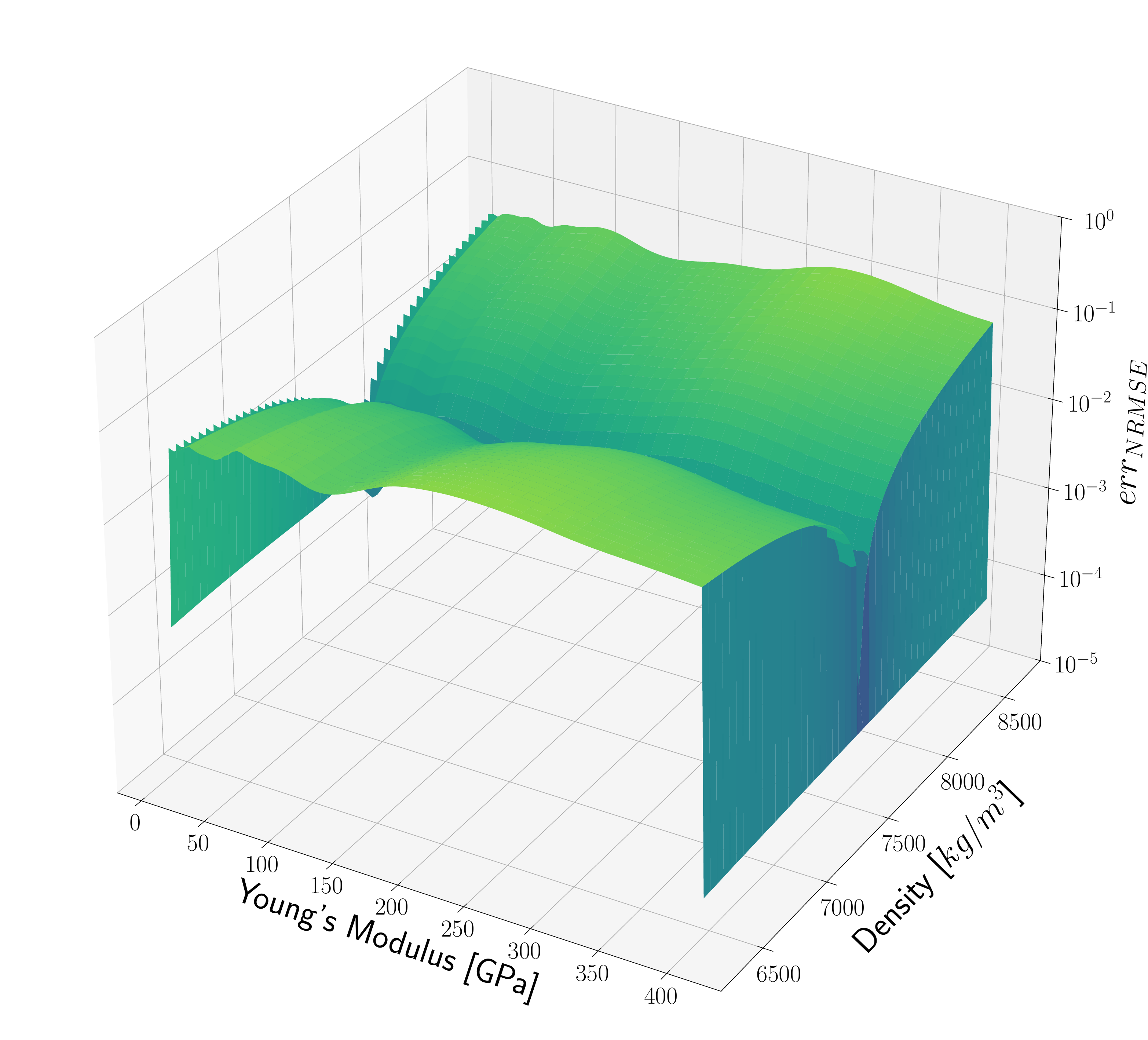}
        \caption{NRMSE for \cite{Panzer2010} , $k=4$, and minMax normalization}
        \label{fig:NRMSE_k4_minMax}
    \end{subfigure}
    \hfill
    \begin{subfigure}[b]{0.3\textwidth}
        \centering
        \includegraphics[width=\textwidth]{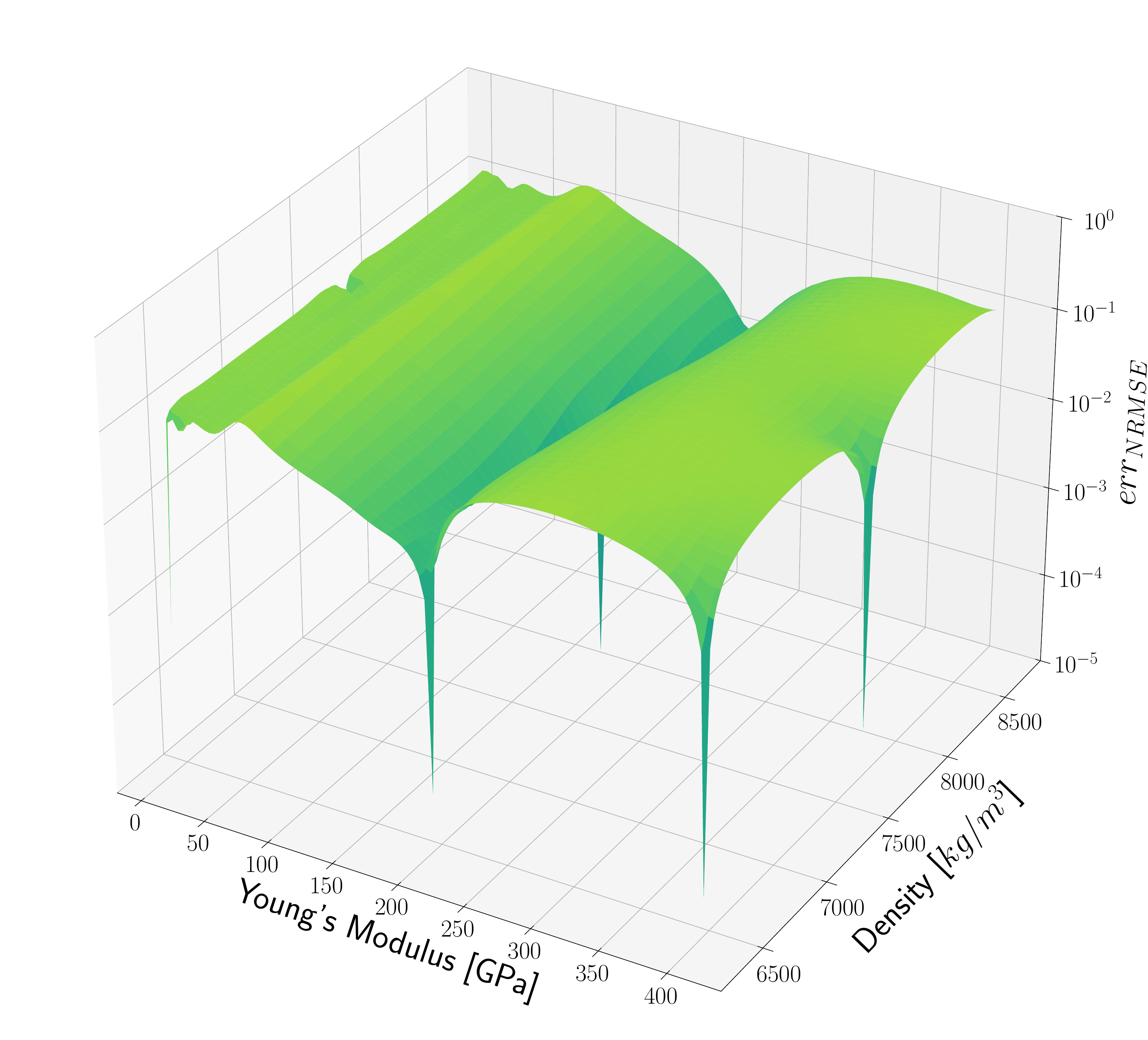}
        \caption{NRMSE for \cite{Panzer2010} , $k=9$, and minMax normalization}
        \label{fig:NRMSE_k9_minMax}
    \end{subfigure}
    \hfill
    \begin{subfigure}[b]{0.3\textwidth}
        \centering
        \includegraphics[width=\textwidth]{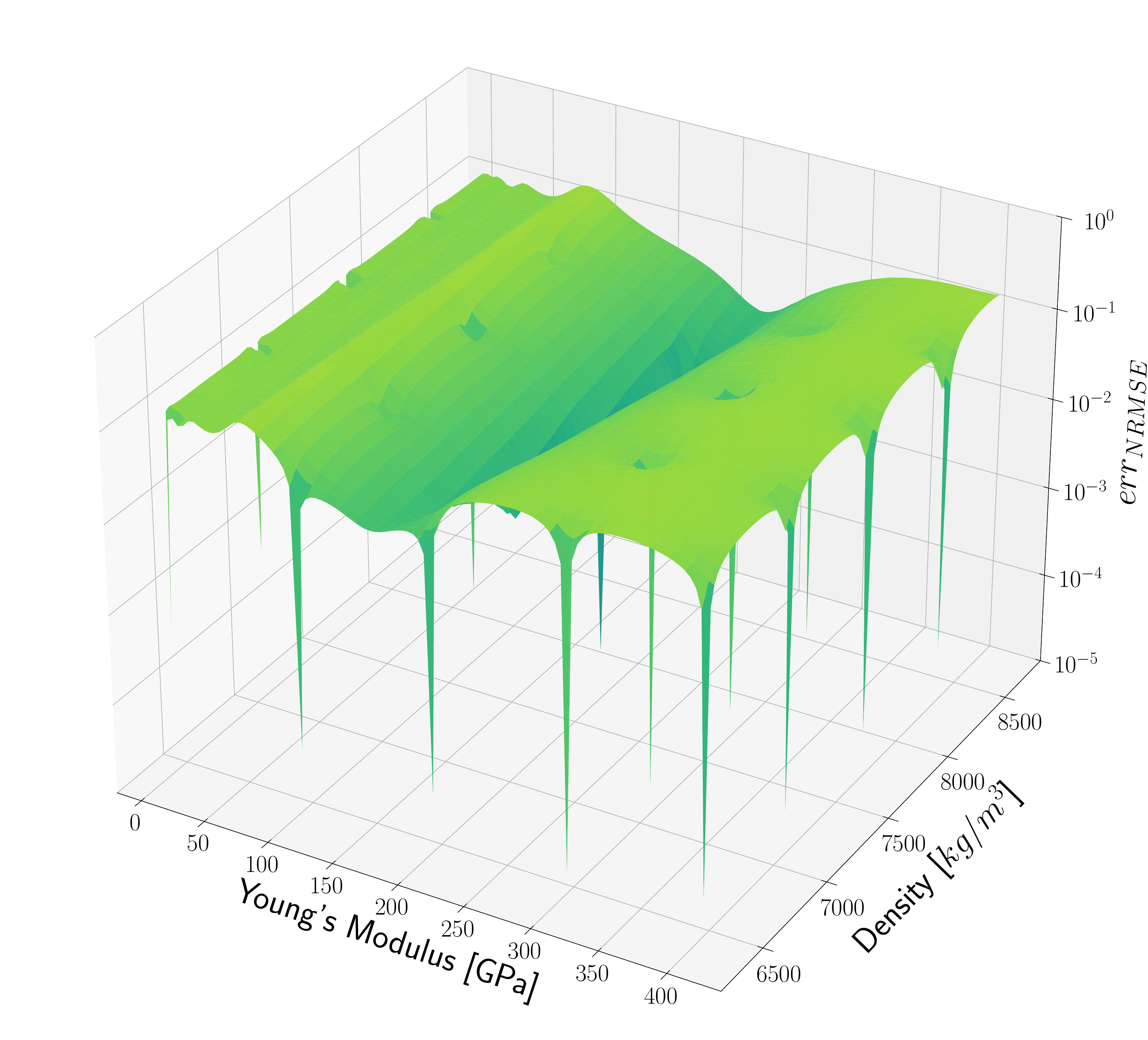}
        \caption{NRMSE for \cite{Panzer2010} , $k=25$, and minMax normalization}
        \label{fig:NRMSE_k25_minMax}
    \end{subfigure}
    \caption{Normalized root mean square error during the fast simulation of the two-parameter surrogate beam model for an increase of training points.}
    \label{fig:SurfaceK}
\end{figure}

\begin{figure}[h]
    \centering 
    \begin{subfigure}[b]{0.28\textwidth}
        \centering
        \includegraphics[width=\textwidth]{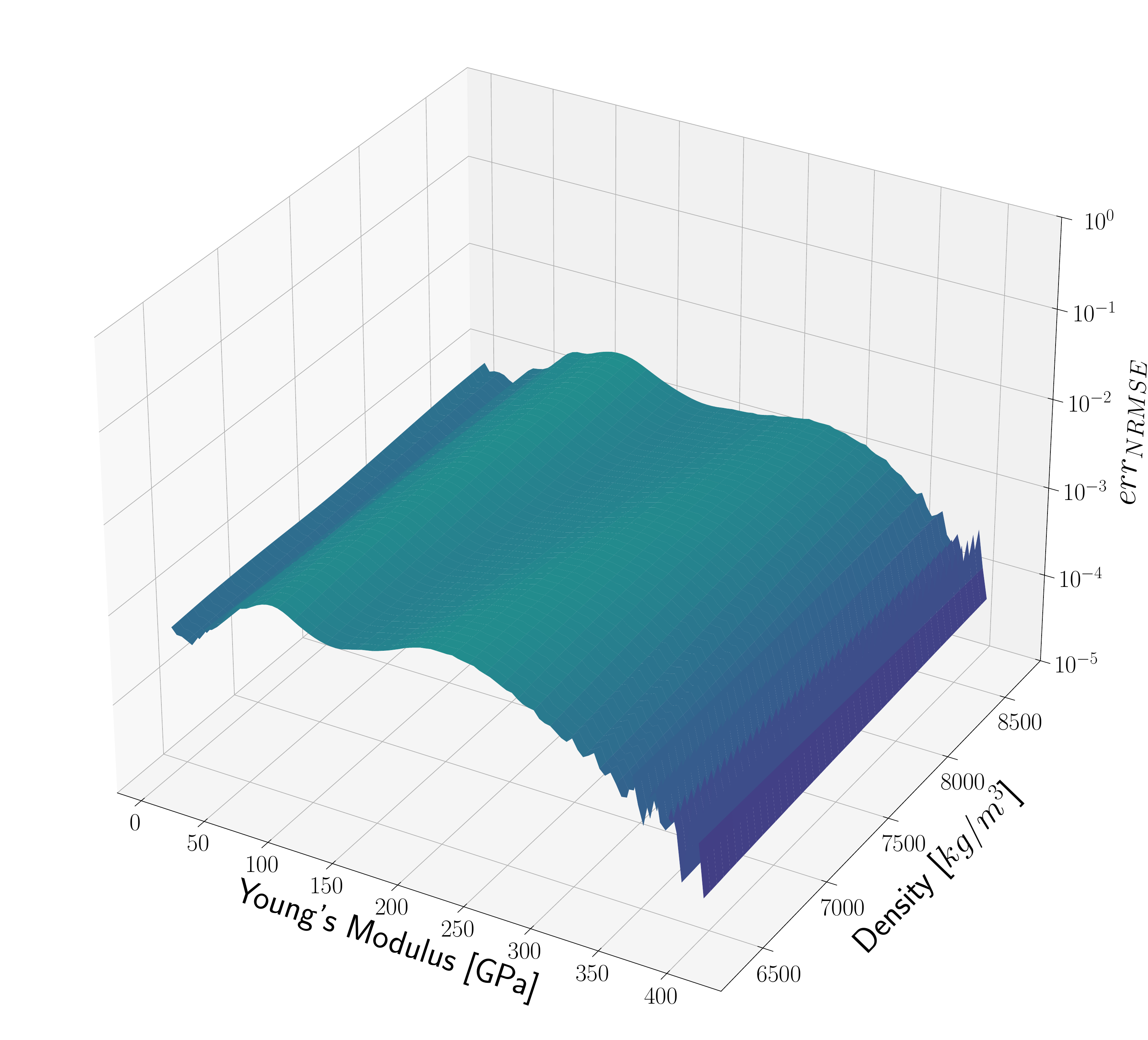}
        \caption{NRMSE for \cite{Panzer2010}, $k=4$ and TPWF}
        \label{fig:NRMSE_k4_TPWF}
    \end{subfigure}
    \hfill
    \begin{subfigure}[b]{0.28\textwidth}
        \centering
        \includegraphics[width=\textwidth]{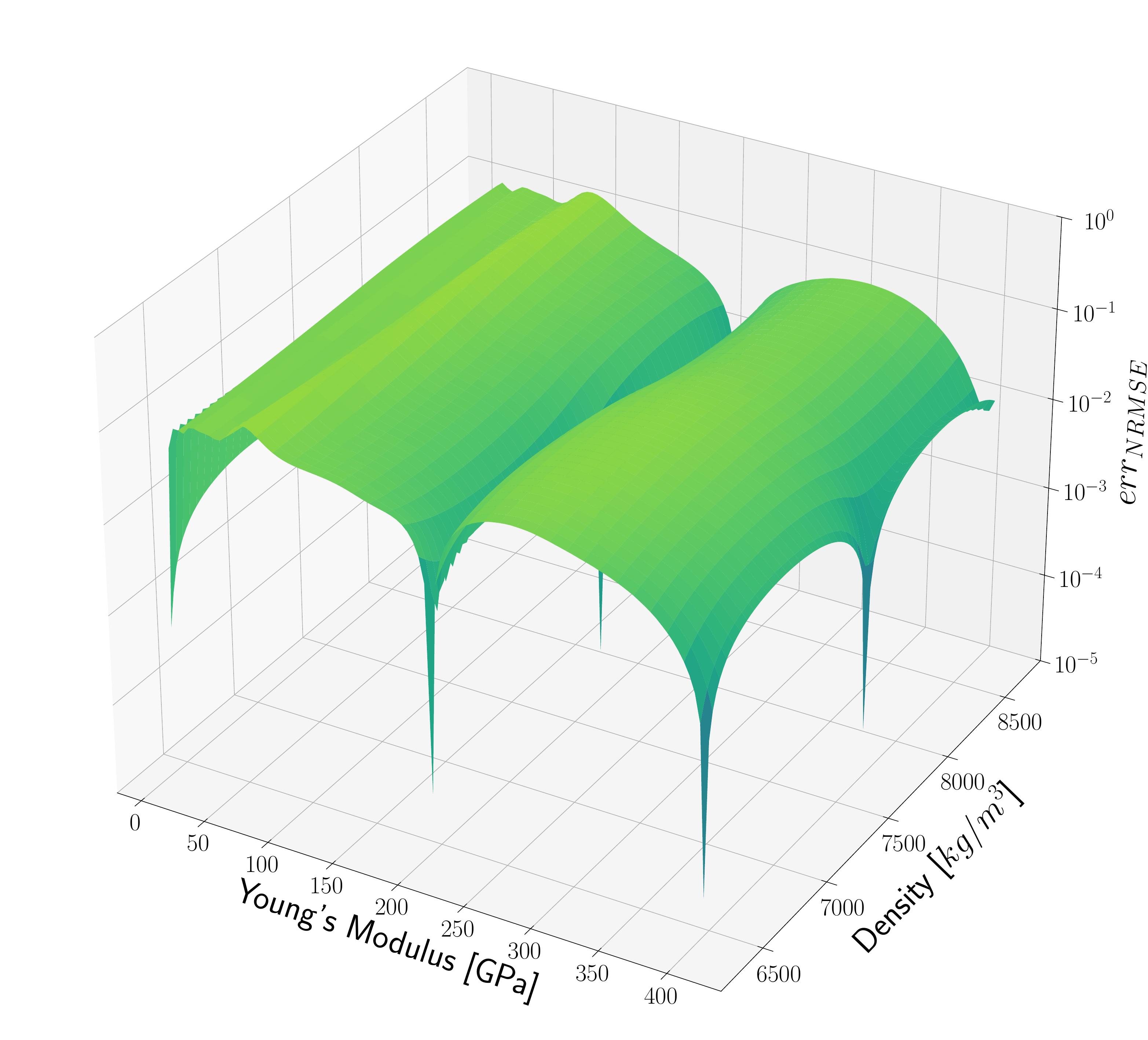}
        \caption{NRMSE for \cite{Panzer2010}, $k=9$ and TPWF}
        \label{fig:NRMSE_k9_TPWF}
    \end{subfigure}
    \hfill
    \begin{subfigure}[b]{0.28\textwidth}
        \centering
        \includegraphics[width=\textwidth]{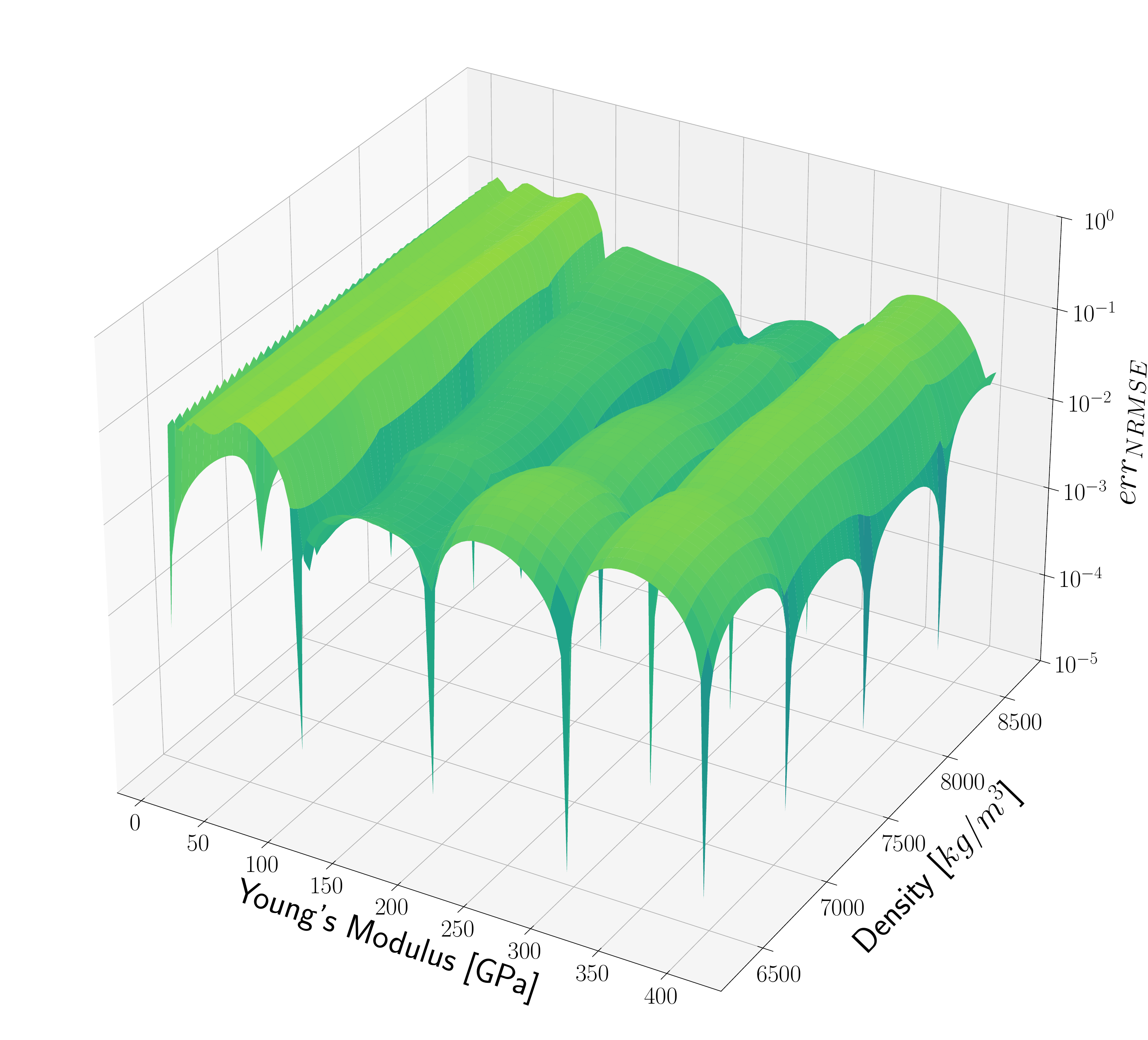}
        \caption{NRMSE for \cite{Panzer2010}, $k=25$ and TPWF}
        \label{fig:NRMSE_k25_TPWF}
    \end{subfigure}
    \caption{Normalized root mean square error surface for tensor product weight function.}
    \label{fig:SurfaceK_TPWF}
\end{figure}

With these parameters, 10201 HFSs are computed to exempplify the presented method by providing a full error surface over the two-dimensional parameter space. The constructed pROMs have $r = 25$ modes. Figure \ref{fig:SurfaceK} depicts the NRMSE over the entire design space, when the method in \cite{Panzer2010} is employed with a minMax normalization and Euclidean distance as weight function for the interpolation for $k=4$, $k=9$, and $k=25$. As expected, a low NRMSE error near the interpolation points is depicted in Fig. \ref{fig:NRMSE_k25_minMax}. However, an increase in the training points does not decrease the error in approximating the displacement field to an acceptable level; the average error for $k=25$ is $14.85\%$. An alternative to this methodology is shown in Fig. \ref{fig:SurfaceK_TPWF}, where the weights for the interpolation are computed using the proposed TPWF, as it eliminates the need for parameter range normalization.

Figure \ref{fig:NRMSE_k4_TPWF} illustrates the NRMSE for $k=4$ and TPWF; the average error is $0.19\%$, an order of magnitude lower than that obtained with the Euclidean weight function. Figures \ref{fig:NRMSE_k9_TPWF} and \ref{fig:NRMSE_k25_TPWF} highlight that an increase of $k$ to nine and 25 does not necessarily lead to improved model accuracy. On the contrary, in most of the parameter space, the accuracy worsens to an average error of $10.93\%$ and $8.18\%$. The increase in training points leads to numerical issues, especially because we are dealing with a large parameter range. The Young's modulus $E$ between training points changes by more than one order of magnitude; subsequently, the entries in the different stiffness matrices $\mathbf{K}_i$ and the reduced stiffness matrices $\mathbf{K}_{r,i}$ differ significantly. Thus, for a $\mu_j$ with a low $E$, as depicted in Fig. \ref{fig:pBI_k9L}, the training points with the highest $E$ will have a bigger influence during the interpolation, leading to an overstiff pROM.

\begin{figure}[ht]
    \centering
    \begin{subfigure}[b]{0.32\textwidth}
        \centering
    \includegraphics[width=\textwidth]{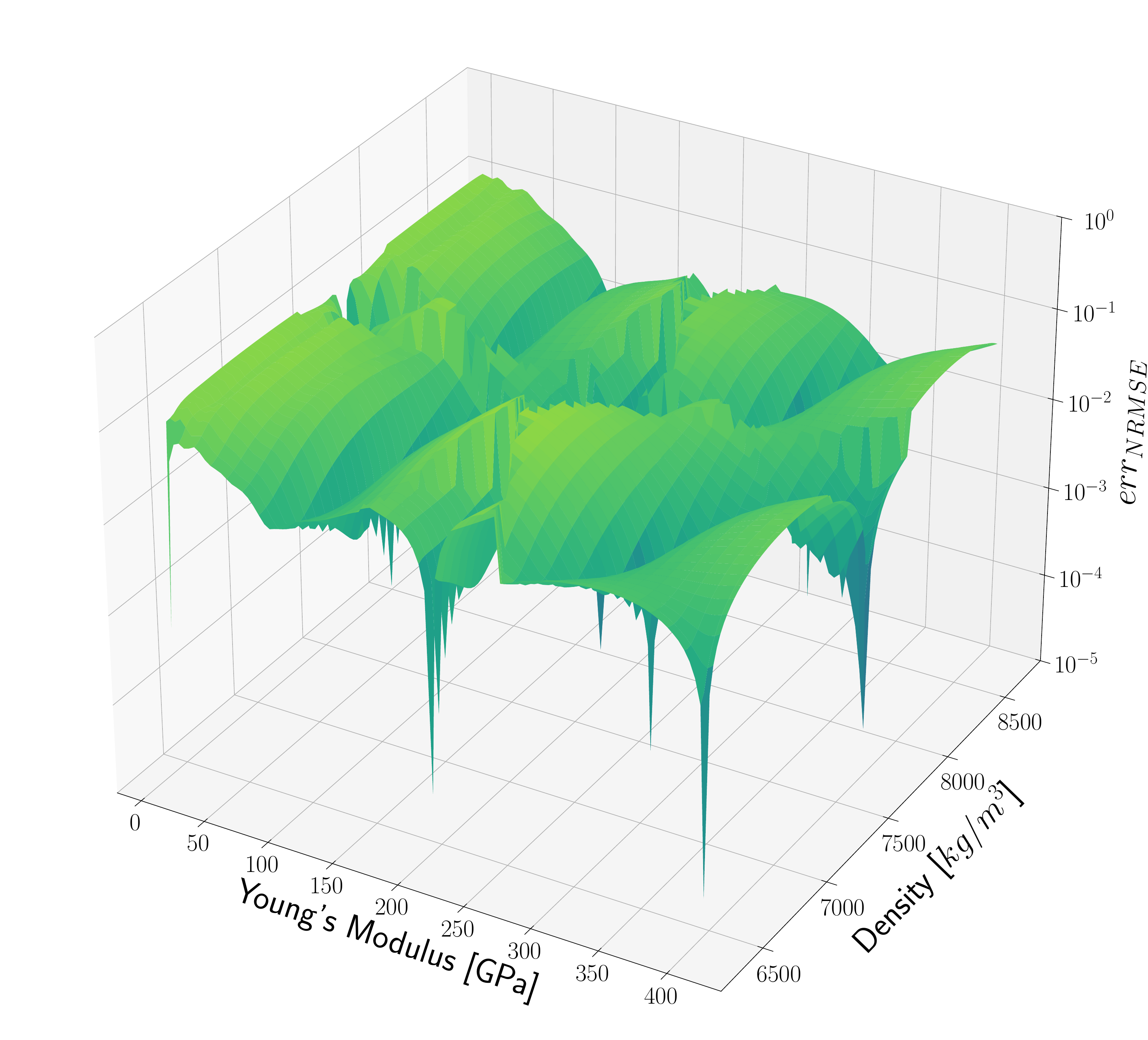}
        \caption{NRMSE for \cite{Yuan2021b}, $k=9$, and Euclidean weight function}
        \label{fig:NRMSE_pBI_k9_kNN_minMax_Euclidean}
    \end{subfigure}
    \hfill
    \begin{subfigure}[b]{0.32\textwidth}
        \centering
        \includegraphics[width=\textwidth]{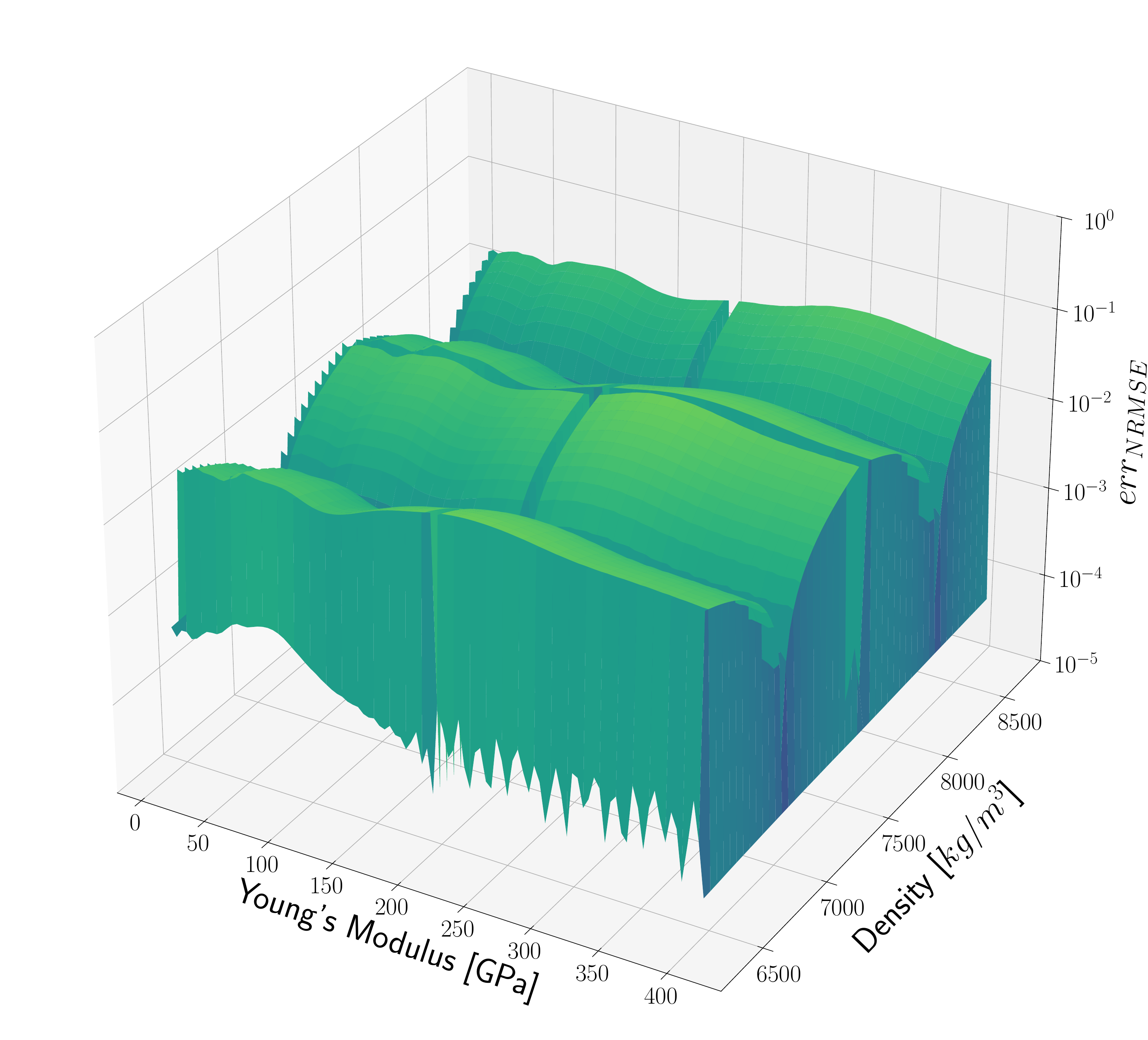}
        \caption{NRMSE for $k=9$, pBI, and Euclidean weight function}
        \label{fig:NRMSE_pBI_k9_Euclidean}
    \end{subfigure}
    \hfill
    \begin{subfigure}[b]{0.32\textwidth}
        \centering
        \includegraphics[width=\textwidth]{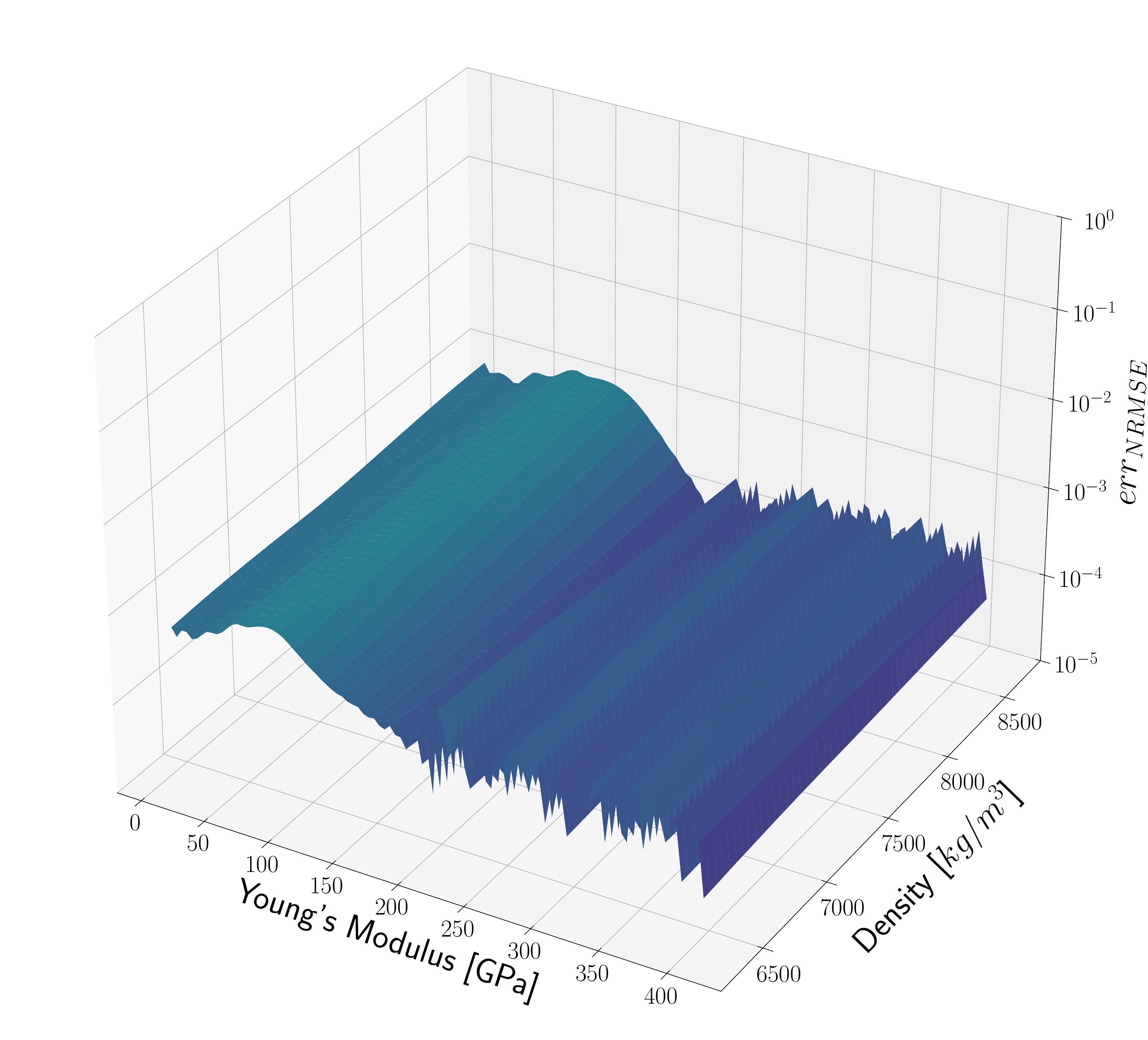}
        \caption{NRMSE for $k=9$, pBI, and TPWF}
        \label{fig:NRMSE_pBI_k9}
    \end{subfigure}
    \caption{Normalized root mean square error surface for different interpolation techniques and weight functions.}
    \label{fig:pBI_Surface}
\end{figure}

A solution to these numerical artifacts for a 1D problem is presented in \cite{Yuan2021, Yuan2021b}, where the interpolation and reduction are performed using only the nearest neighbors, and a Euclidean function computes the weights $\omega_i$. Figure \ref{fig:NRMSE_pBI_k9_kNN_minMax_Euclidean} illustrates the surface error of the technique proposed in \cite{Yuan2021b}, when applied only during the interpolation stage, showing a lower error than the method presented in \cite{Panzer2010} (Fig. \ref{fig:SurfaceK}). Similarly, Figure \ref{fig:NRMSE_pBI_k9_Euclidean} presents the surface error for the method introduced in Sec. \ref{ss:BI}, which employs the box clustering algorithm and the Euclidean function for the weights $\omega$. However, both methods have a larger NRMSE than the one in Fig. \ref{fig:NRMSE_k4_TPWF}. Figure \ref{fig:NRMSE_pBI_k9} displays the NRMSE of pBI combined with TPWF, showing that pBI reduces the error across the entire parameter space. By not computing the pROM with all the training points, the artificial stiffness due to the influence of distant training points is reduced. The average error of the new method, pBI with TPWF, is $0.056\%$. Figure \ref{fig:pBI_Surface} demonstrates that TPWF has a better performance than conventional approaches, which include normalization. 

\begin{figure}[ht]
    \centering
    \begin{subfigure}[b]{0.32\textwidth}
        \centering
        \includegraphics[width=\textwidth]{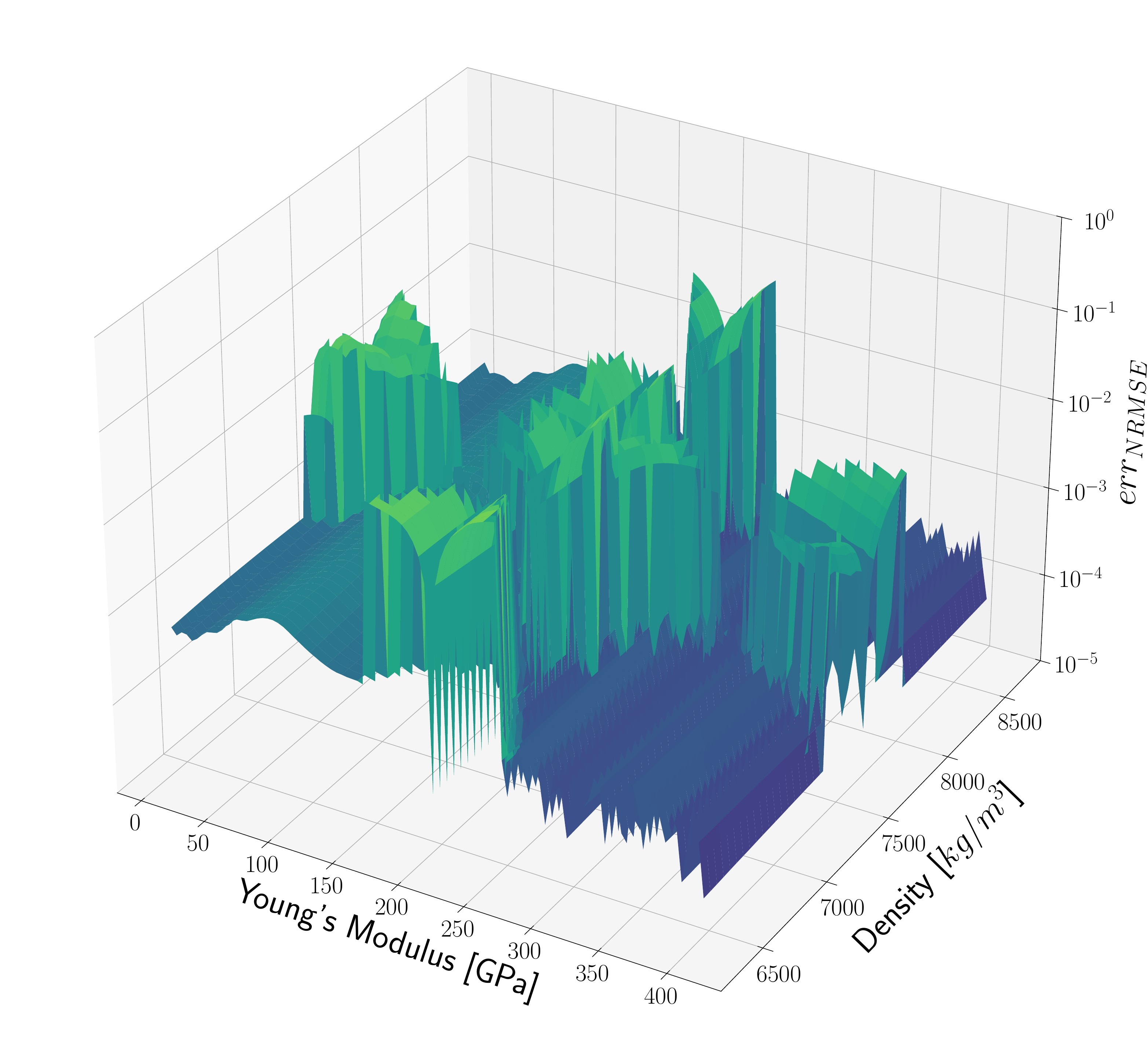}
        \caption{NRMSE for \cite{Yuan2021b}, $k=9$, and kNN}
        \label{fig:NRMSE_pBR_k9_kNN}
    \end{subfigure}
    \hfill
    \begin{subfigure}[b]{0.32\textwidth}
        \centering
        \includegraphics[width=\textwidth]{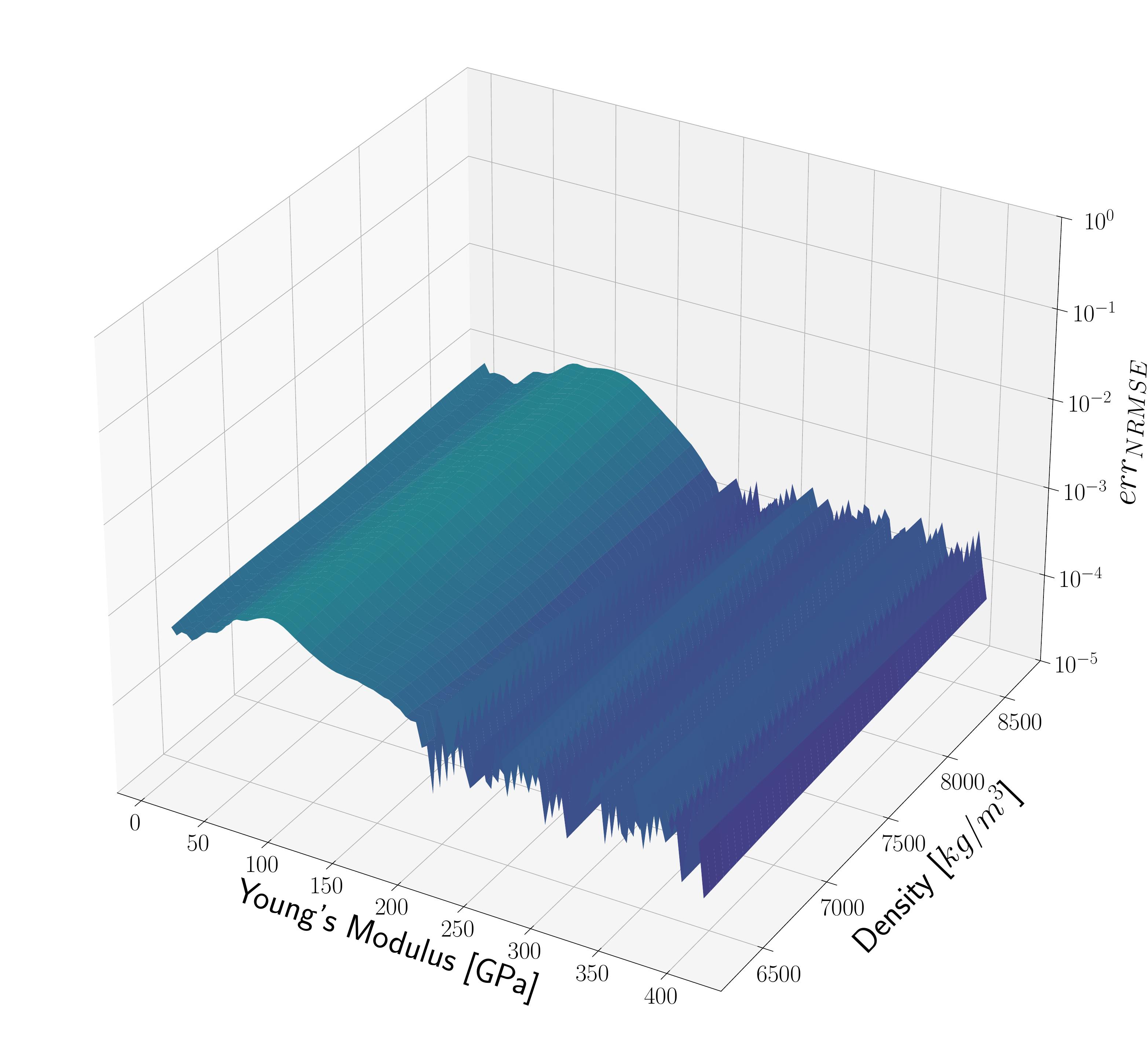}
        \caption{NRMSE for $k=9$ and pBR}
        \label{fig:NRMSE_pBR_k9}
    \end{subfigure}
    \hfill
    \begin{subfigure}[b]{0.32\textwidth}
        \centering
        \includegraphics[width=\textwidth]{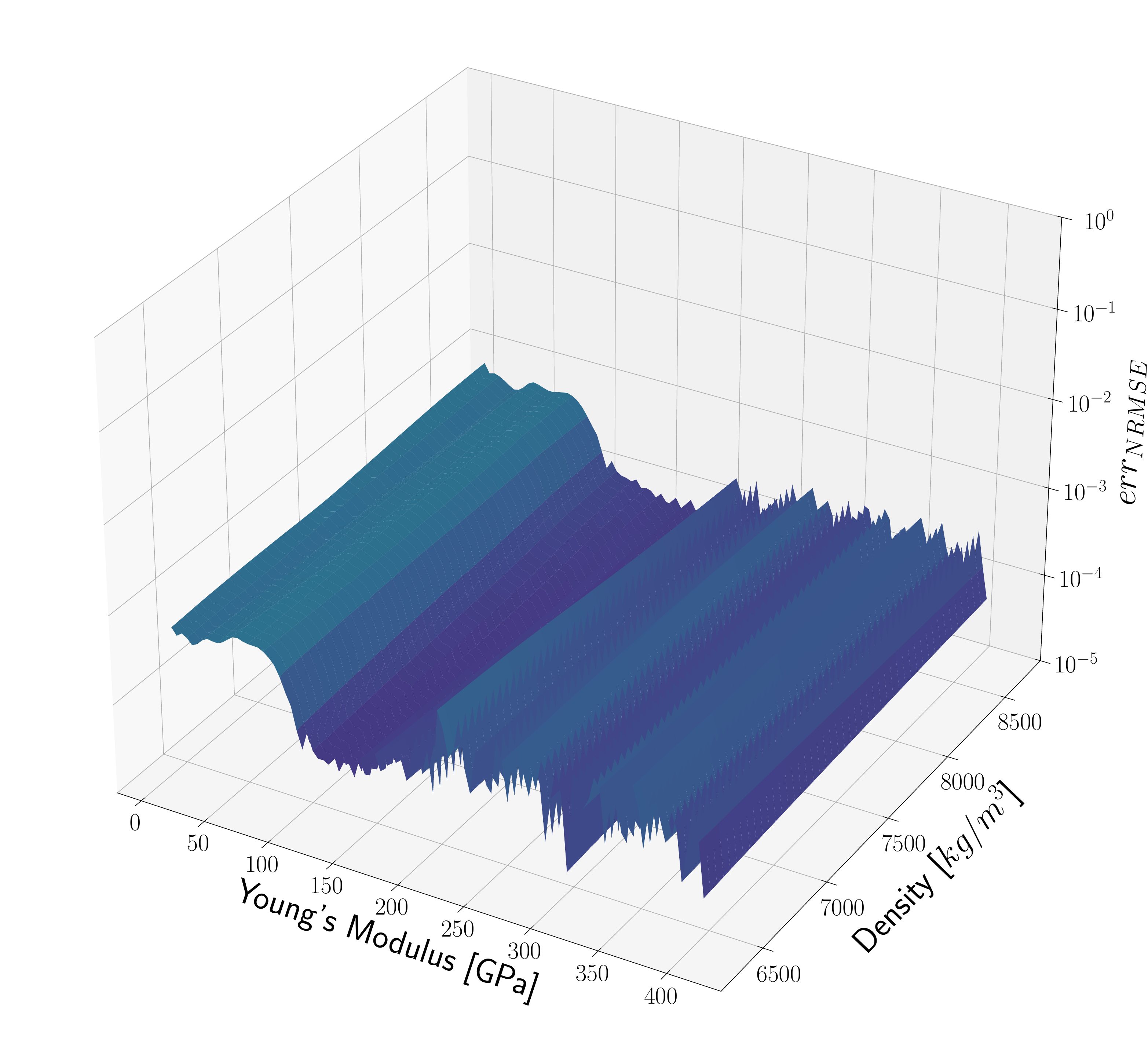}
        \caption{NRMSE for $k=25$ and pBR}
        \label{fig:NRMSE_pBR_k25}
    \end{subfigure}
    \caption{Normalized root mean square error during the fast simulation of the two-parameter surrogate beam model employing parametric box reduction with kNN and box clustering.}
    \label{fig:pBR_Surface}
\end{figure}

\begin{figure}[ht]
    \centering
    \begin{subfigure}[b]{0.24\textwidth}
        \centering
        \begin{tikzpicture}
            \node[anchor=south west,inner sep=0] (image) at (0,0) {\includegraphics[width=\textwidth]{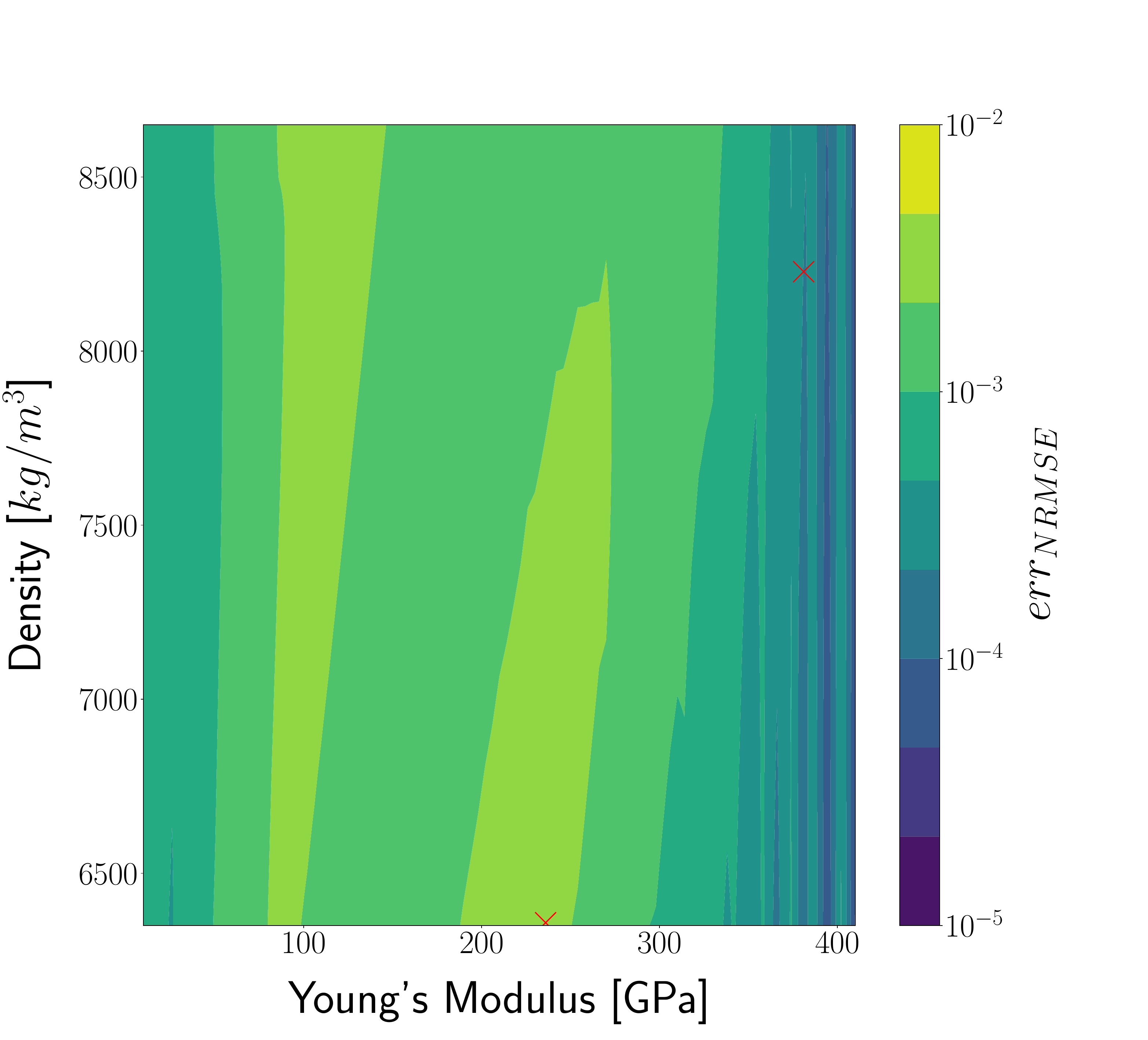}};
            \begin{scope}[x={(image.south east)},y={(image.north west)}]
            \node[anchor=west, black] at (0.6,0.7) {\tiny $\mu_{2}$};
            \node[anchor=west, black] at (0.4,0.15) {\tiny $\mu_{1}$};
            \end{scope}
         \end{tikzpicture}
        \caption{NRMSE for \cite{Panzer2010}, $k=4$, and TPWF}
        \label{fig:NRMSE_k4_TPWF_contour}
    \end{subfigure}
    \hfill
    \begin{subfigure}[b]{0.24\textwidth}
        \centering
        \begin{tikzpicture}
            \node[anchor=south west,inner sep=0] (image) at (0,0) {\includegraphics[width=\textwidth]{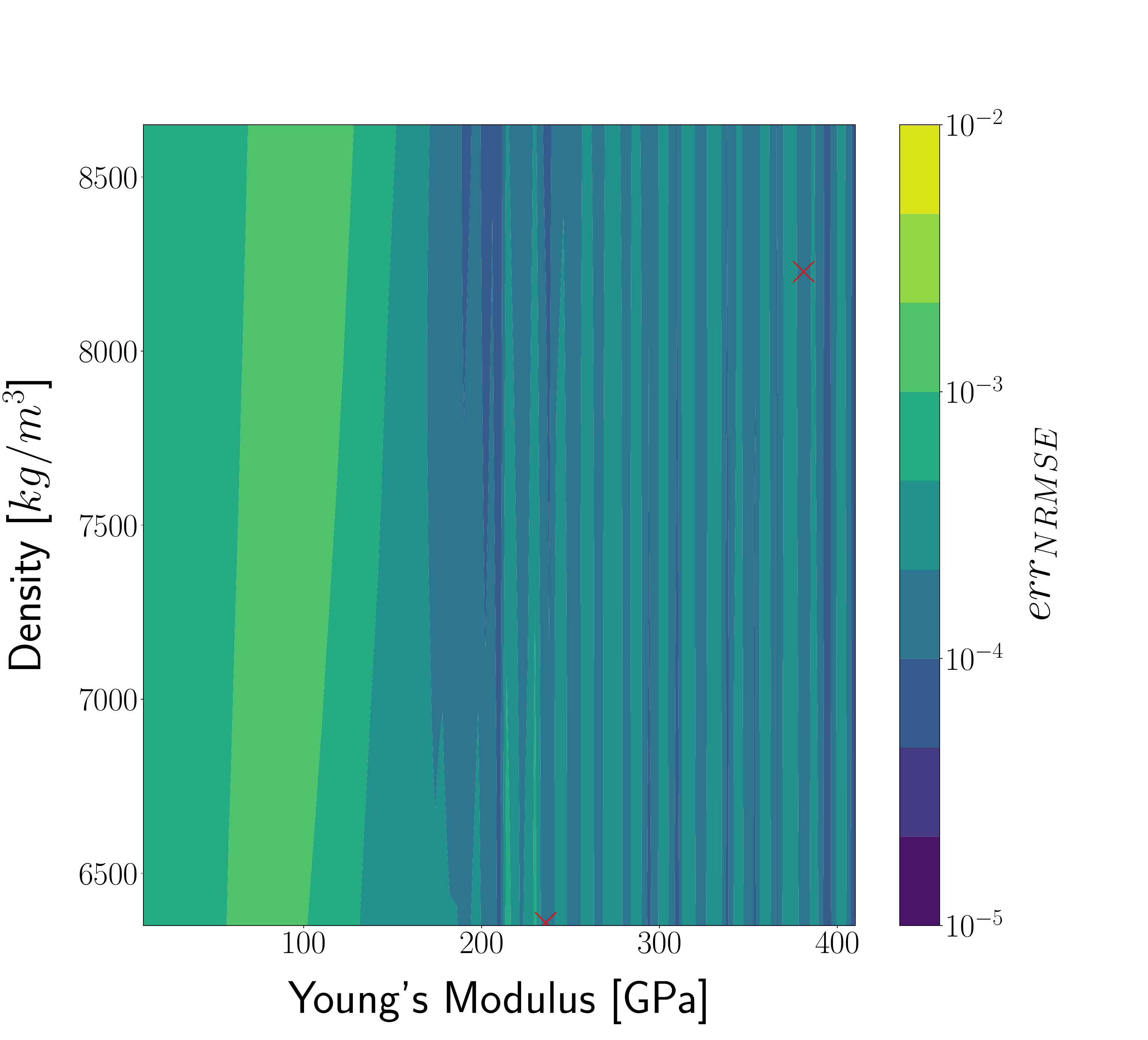}};
            \begin{scope}[x={(image.south east)},y={(image.north west)}]
            \node[anchor=west, black] at (0.6,0.7) {\tiny $\mu_{2}$};
            \node[anchor=west, black] at (0.4,0.15) {\tiny $\mu_{1}$};
            \end{scope}
         \end{tikzpicture}
        
        \caption{NRMSE for pBI, $k=9$, and TPWF}
        \label{fig:NRMSE_pBI_k9_contour}
    \end{subfigure}
    \hfill
    \begin{subfigure}[b]{0.24\textwidth}
        \centering
        \begin{tikzpicture}
            \node[anchor=south west,inner sep=0] (image) at (0,0) {\includegraphics[width=\textwidth]{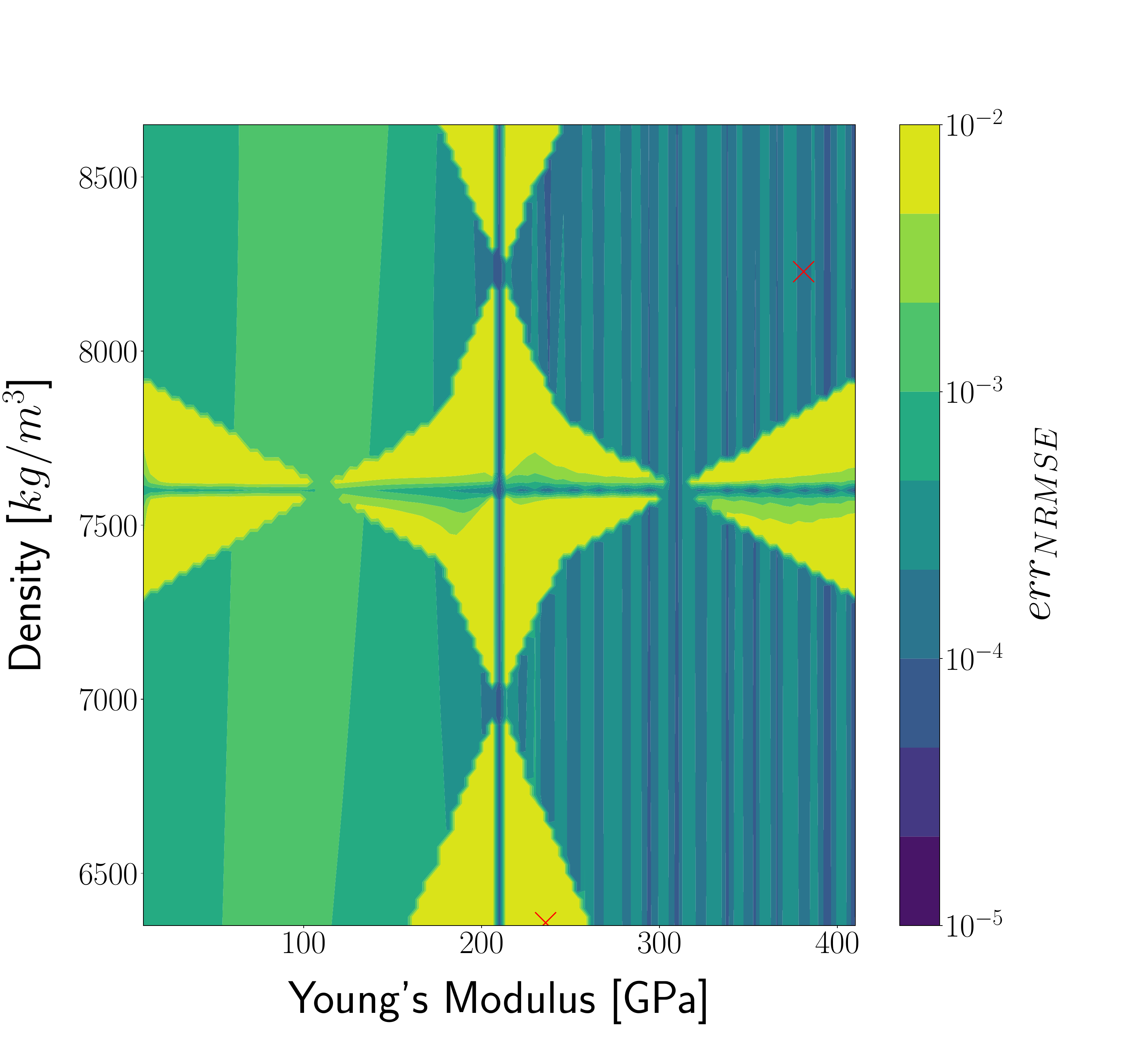}};
            \begin{scope}[x={(image.south east)},y={(image.north west)}]
            \node[anchor=west, black] at (0.6,0.7) {\tiny $\mu_{2}$};
            \node[anchor=west, black] at (0.4,0.15) {\tiny $\mu_{1}$};
            \end{scope}
         \end{tikzpicture}
        
        \caption{NRMSE for pBR, $k=9$, and pBR}
        \label{fig:NRMSE_pBR_k9_kNN_contour}
    \end{subfigure}
    \hfill
    \begin{subfigure}[b]{0.24\textwidth}
        \centering
        \begin{tikzpicture}
            \node[anchor=south west,inner sep=0] (image) at (0,0) {\includegraphics[width=\textwidth]{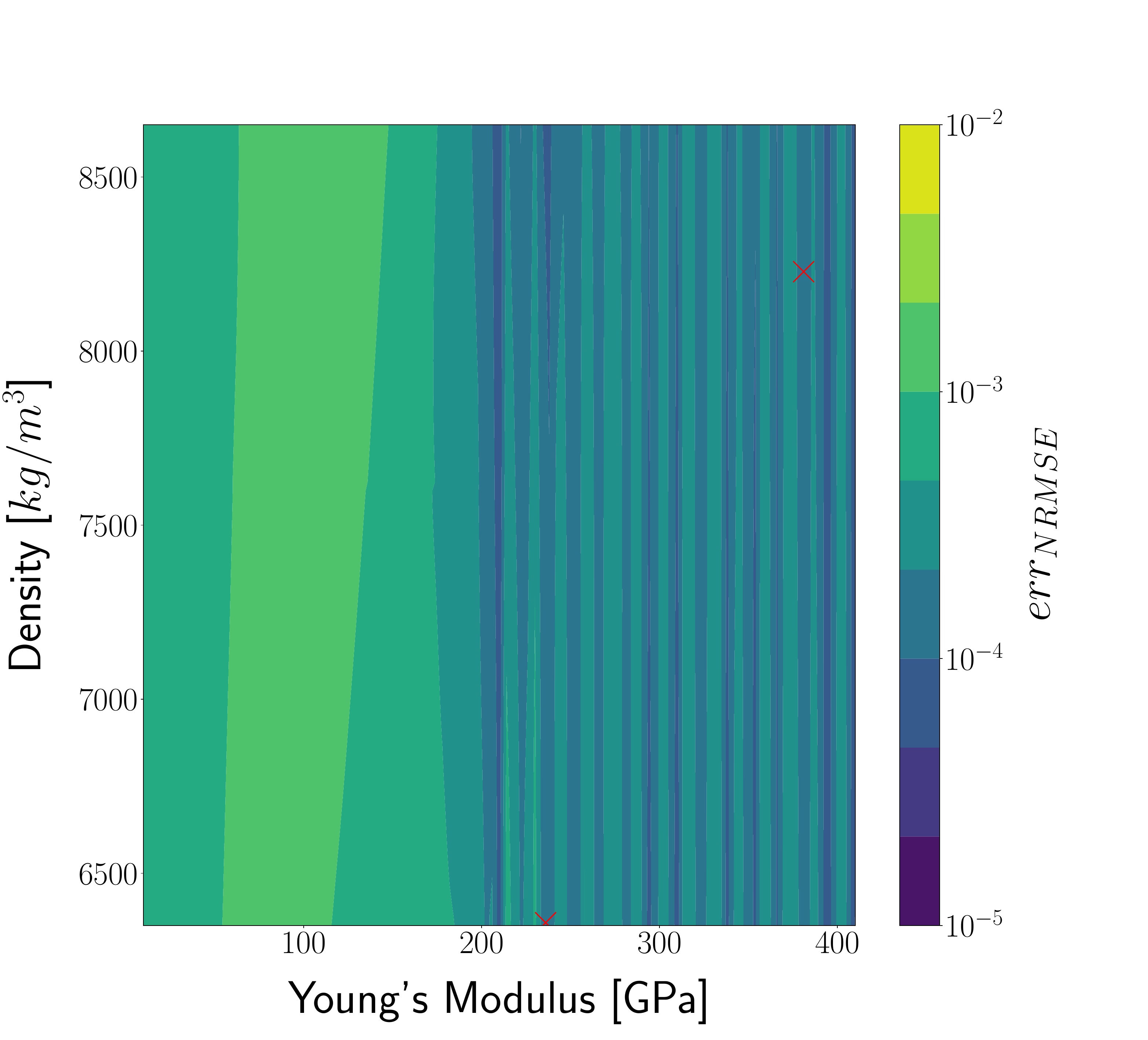}};
            \begin{scope}[x={(image.south east)},y={(image.north west)}]
            \node[anchor=west, black] at (0.6,0.7) {\tiny $\mu_{2}$};
            \node[anchor=west, black] at (0.4,0.15) {\tiny $\mu_{1}$};
            \end{scope}
         \end{tikzpicture}
        
        \caption{NRMSE for $k=9$ and pBR}
        \label{fig:NRMSE_pBR_k9_contour}
    \end{subfigure}
    \caption{Normalized root mean square error during the fast simulation of the two-parameter surrogate beam model employing parametric box reduction with kNN and box clustering.}
    \label{fig:pBR_contour}
\end{figure}

Similarly, pBR can be employed to remedy the overstiff behavior. Figure \ref{fig:pBR_Surface} depicts the surface error of the main method of this work, and how it is contrasted when kNN is applied instead of box clustering. Given that TPWF showcases a better error performance (Figs. \ref{fig:SurfaceK_TPWF} and \ref{fig:pBI_Surface}), the weights $\omega$ are computed with TPWF. Additionally, the contour plots of the NRMSE associated with the best-performing methods are presented in Fig. \ref{fig:pBR_contour}. The method proposed in \cite{Yuan2021b} with TPWF exhibits behavior opposite to that of normal pROM techniques. It has a higher error near the interpolation points. This phenomenon is depicted in Figs. \ref{fig:NRMSE_pBR_k9_kNN} and \ref{fig:NRMSE_pBR_k9_kNN_contour}; it is one of the downsides of using kNN as a clustering algorithm, as explained in Fig. \ref{fig:NearPoints}. Similarly, pBR can be employed to remedy the overstiff behavior. As highlighted Figs. \ref{fig:NRMSE_pBR_k9} and \ref{fig:NRMSE_pBR_k9_contour}, both methods achieve the same accuracy level for this particular test case. 

\begin{figure}[ht]
     \centering
     \begin{subfigure}[b]{0.48\textwidth}
         \centering
         \includegraphics[width=\textwidth]{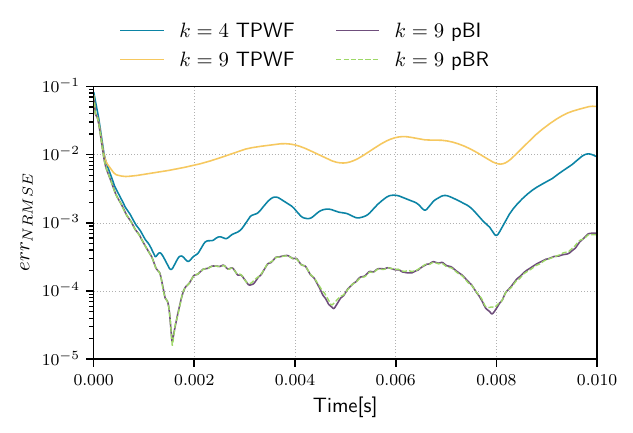}
         \caption{NRMSE for $\mu_{1}$ for classical approach and pBI for four a nine training points}
         \label{fig:B_p2_MSRE_1}
     \end{subfigure}
     \hfill
     \begin{subfigure}[b]{0.48\textwidth}
         \centering
         \includegraphics[width=\textwidth]{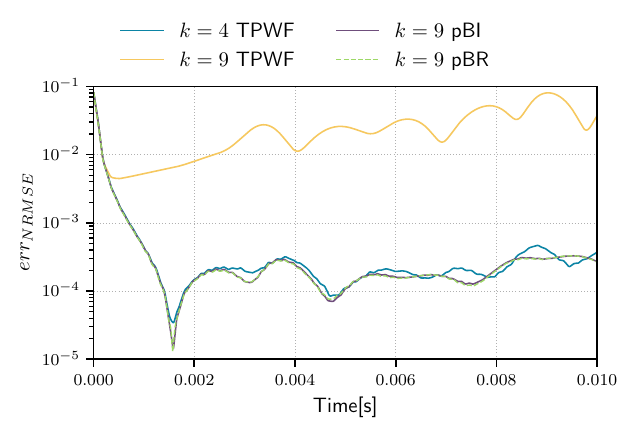}
         \caption{NRMSE for $\mu_{2}$ for classical approach and pBI for four a nine training points}
         \label{fig:B_p2_MSRE_2}
     \end{subfigure}
        \caption{Comparison of accuracy across different methods by analyzing the evolution of normalized root mean square error through the simulation time at two separate validation points, $\mu_{1}=[236\, \mathrm{GPa}, 6358\, \frac{\mathrm{kg}}{\mathrm{m^3}}]$ and $\mu_{2}=[381\, \mathrm{GPa}, 8228\, \frac{\mathrm{kg}}{\mathrm{m^3}}]$.}
        \label{fig:B_p2_nrmse}
\end{figure}

Furthermore, figure \ref{fig:B_p2_nrmse} analyzes the reduced models' behavior for two validation points, $\mu_{1}$ and $\mu_{2}$. The x-axis denotes the simulation time progression, whereas the y-axis measures the magnitude of the NRMSE. For $\mu_{1}=[236\, \mathrm{GPa}, 6358\, \frac{\mathrm{kg}}{\mathrm{m^3}}]$. In Fig. \ref{fig:B_p2_MSRE_1}, the accuracy decreases when more training points are selected. The NRMSE is $0.93\%$ and $5.1\%$ for four and nine training points, respectively. Additionally, the pBI and pBR exhibit the best accuracy, with an NRMSE of $0.07\%$. To understand the phenomena behind these results, Fig. \ref{fig:pBI_k4L} and \ref{fig:pBI_k9L} illustrate the parameter space for four and nine training points. For nine training points, $\mu_2$ is expected to have the biggest influence on the reduced model. However, given the large parameter ranges, the entries of the system matrices of training point $\mu_9$ are much larger than the entries from $\mu_2$. Therefore $ \omega_2\textbf{X}_{r,2} \ll \omega_9 \textbf{X}_{r,9}$, even though $\omega_2 > \omega_9$. Furthermore, pBI and pBR allow an increase in accuracy by having more training points, without the numerical artifacts. Figure \ref{fig:NRMSE_pBR_k25} depicts the surface error for $k=25$, it achieves an average NRMSE of $0.035\%$ over the parameter space, corresponding to a $35\%$ decrease compared to $k=9$. Similarly, Fig. \ref{fig:B_p2_MSRE_2} illustrates the behavior for the validation point $\mu_{2}=[381\, \mathrm{GPa}, 8228\, \frac{\mathrm{kg}}{\mathrm{m^3}}]$. It is shown that the classical approach with nine training points delivers the highest value of error measure, as other training points create a non-physical impact on the interpolation. In contrast to the previous validation point, pBI and pBR do not improve the accuracy compared to using only four training points. The main reason for this behavior is that $\mu_{2}$ is rather close to a training point with high parameter values (Fig. \ref{fig:pBI_k4H}).

\begin{figure}[ht]
    \centering
    \begin{tikzpicture}
        \node[anchor=south west,inner sep=0] (image) at (0,0) {\includegraphics[width=0.5\linewidth]{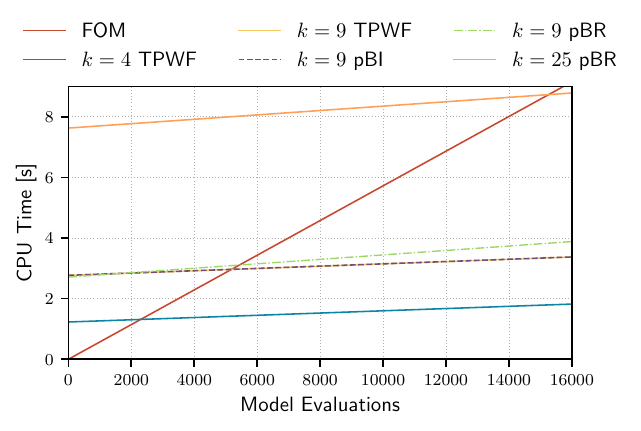}};
        \begin{scope}[x={(image.south east)},y={(image.north west)}]  
            \tikzstyle{arrow} = [->,>=stealth]
            \draw[arrow] (0.45,0.6) -- (0.8,0.76);
            \node[black] at (0.4,0.55) {\shortstack{\tiny\textbf{Break-even} \\ \tiny \textbf{point}}};
            \draw[arrow] (0.41,0.52) -- (0.38,0.4);
        \end{scope}
    \end{tikzpicture}
    \caption{Computational time comparison between the FOM and the different reduction techniques, the initial value of the pROM corresponds to the offline phase time}
    \label{fig:B_time}
\end{figure}

The previous results showed a detailed evaluation of the accuracy of the pBI and pBR against methods presented in \cite{Panzer2010, Yuan2021, Yuan2021b}. To complement the analysis, figure \ref{fig:B_time} presents the trade-off between the accuracy and computational effort. The horizontal axis represents the number of model evaluations (i.e., time integration steps), while the vertical axis corresponds to computational time. The required time for the model creation (i.e., offline phase) is depicted in the plot as the CPU time for 0 model evaluations. Additionally, the break-even point is indicated on the plot. For four training points and computing the interpolation matrices with TPWF the break-even point is $2297$; while for nine training points, it is $5182$ for the method in \cite{Panzer2010}, $5184$ for pBI, and $5404$ for pBR. Therefore, no additional computational burden is added with the newly introduced method compared to classical pROM approaches. The decrease in computational efficiency by increasing the number of training points is justified by the gain in accuracy, as shown in Fig. \ref{fig:pBR_contour}. Given that the system reduction time grows linearly with the number of training points. Whereas the basis change time increases non-linearly because the SVD performed on $V_{all}$. Therefore, the required CPU time for the offline stage is significantly affected by the number of training points. Furthermore, the pBI CPU time associated with solving the reduced system and projecting back the solution to the full order space (i.e., online phase) is reduced by a factor of 15 compared to the HFS solving time, corresponding to a speed-up of 15. Even though the speed-up for pBR is in the same range as for pBI, the CPU time for pBR is higher than that of pBI because of the required computation of an additional SVD at each new parameter. 

\subsection{Power module}

\begin{figure}[h]
    \centering
    \begin{subfigure}[b]{0.48\textwidth}
        \centering
        \begin{tikzpicture}
            \node[anchor=south west,inner sep=0] (image) at (0,0) {\includegraphics[width=0.95\linewidth]{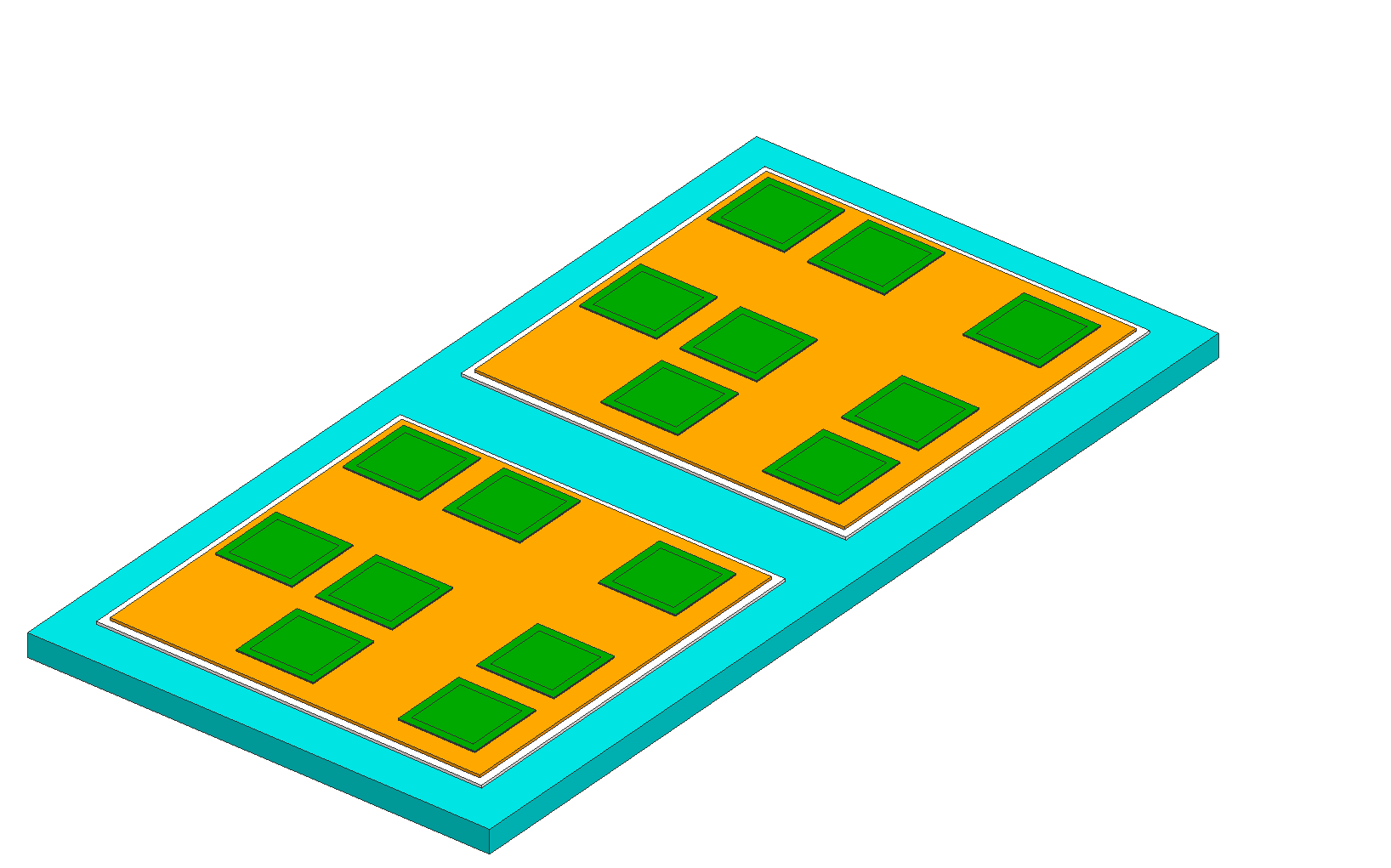}};
            \begin{scope}[x={(image.south east)},y={(image.north west)}]  
                \draw[black,thick,->] (0.86,0.75) -- (0.74,0.62);
                \node[anchor=west, black] at (0.85,0.75) {\tiny Body 1};
                \draw[black,thick,->] (0.8,0.35) -- (0.52,0.28);
                \draw[black,thick,->] (0.8,0.35) -- (0.63,0.4);
                \node[anchor=west, black] at (0.8,0.35) {\tiny Body 2};
            \end{scope}
        \end{tikzpicture}
        \caption{Power module 3D model}
        \label{fig:PowerModule}
    \end{subfigure}
    \begin{subfigure}[b]{0.48\textwidth}
        \centering
        \includegraphics[width=0.95\linewidth]{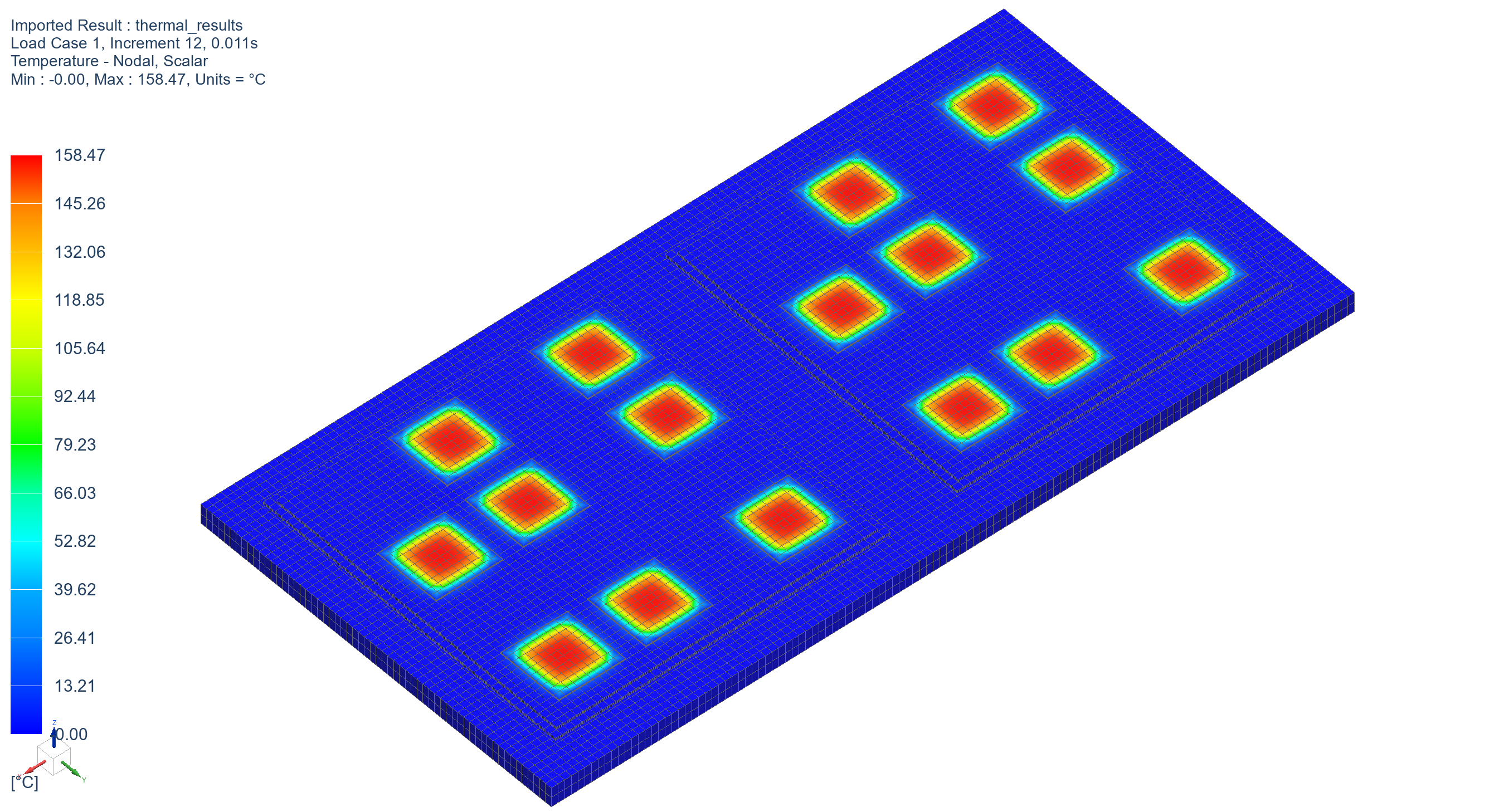}
        \caption{Temperature field HFS at $t=0.011\; \mathrm{s}$}
        \label{fig:HFS}
    \end{subfigure}
    \caption{3D power module CAD model and high fidelity temperature field snapshot computed with Simcenter 3D thermal solver}
    \label{fig:Sim3DPM}
\end{figure}

Power modules (Fig. \ref{fig:PowerModule}) are a core component of a power converter; they serve as a physical container for several semiconductor switches, typically including diodes and Insulated Gate Bipolar Transistor (IGBT) dies. During the operation, the power modules generate heat and must be cooled to avoid damage to the components and maintain the temperature within the optimal operational range \cite{Kulkarni2023}. During the design phase, multiple power modules are created by varying some geometrical (thickness, position) and material parameters (thermal conductivity, mass density, specific heat).

For this purpose, the transient thermal behavior of the design needs to be evaluated to select the optimal design. However, each HFS is computationally expensive, creating a bottleneck during the optimization. Therefore, building a pROM is of great interest for this industrial application case, as the design engineer can quickly evaluate the transient thermal behavior of different power module designs. Figure \ref{fig:PowerModule} displays a power module with 16 chips (Body 1) and two substrates. A heat load to a section of the chips and a convection boundary condition to the baseplate's bottom face are applied. All the thermal HFSs are performed using Simcenter 3D and its thermal solver \cite{Reiner2022}. From the HFS, we extract the system matrices and vectors to perform the system reduction. With this information, a first-order LTI system is assembled as shown in Eq. \eqref{eq:FOM_1st}. The power module FOM has $n = 244,577$ DoFs, and the FOM solving time is $t_{FOM}=2.602\;\mathrm{s}$ for each implicit Euler iteration. Figure \ref{fig:HFS} illustrates the temperature field at a $t=0.011\; \mathrm{s}$ obtained with Simcenter 3D.

To evaluate the accuracy of the model over the entire domain, we calculate the Mean Square Relative Error as an indicator

\begin{equation}\label{eq:MSRE}
    MSRE = \sqrt{\frac{\sum\limits_{i=1}^{n}{\left(\frac{u_s- u_{r,i}}{u_s}\right)^2}}{n}} \,.
\end{equation}

In this experiment, a two-parameter surrogate model is created for a power module to test the method's applicability in industrial scenarios. Similarly, the power module's FOM dimension is 407 times larger than that of the beam model, serving as a test case to evaluate the reduction technique's capability to handle significantly larger systems and assess its versatility with first-order systems. The surrogate's design variables are the thermal conductivity of body 1 and body 2 (Fig. \ref{fig:PowerModule}) and their ranges are $K_1 = [100, 300]\frac{\mathrm{W}}{\mathrm{mK}}$ and $K_2 = [20, 290]\frac{\mathrm{W}}{\mathrm{mK}}$. For validation purposes, 29 designs are sampled with a 2D LHS, and for each, an HFS is carried out to compare their results with the pROMs. 

Based on the previous experiment, pBI generates more accurate models than the classical approach when $k>2^p$, i.e., additional points for the training set. Thus, for the two-parameter power module model, we use pBI when more than four training points are available, the minimum $k$ required to create the structured grid. For $k=4$, all the methods are equivalent. Furthermore, for investigating the influence of additional training points, a pROM is computed with a training set sampled using a structured $3\times4$ grid. A denser discretization of $K_2$ in the structured grid is selected, given that the parameter range for $K_2$ is larger than for $K_1$. The reduced models have $100$ DoFs; this value is chosen since higher values do not lead to significant model improvement, and the pROM dimension maintains a low-dimensional representation for fast calculations. 

\begin{figure}[ht]
     \centering
     \begin{subfigure}[b]{0.46\textwidth}
         \centering
         \begin{tikzpicture}
            \node[anchor=south west,inner sep=0] (image) at (0,0) {\includegraphics[width=\textwidth]{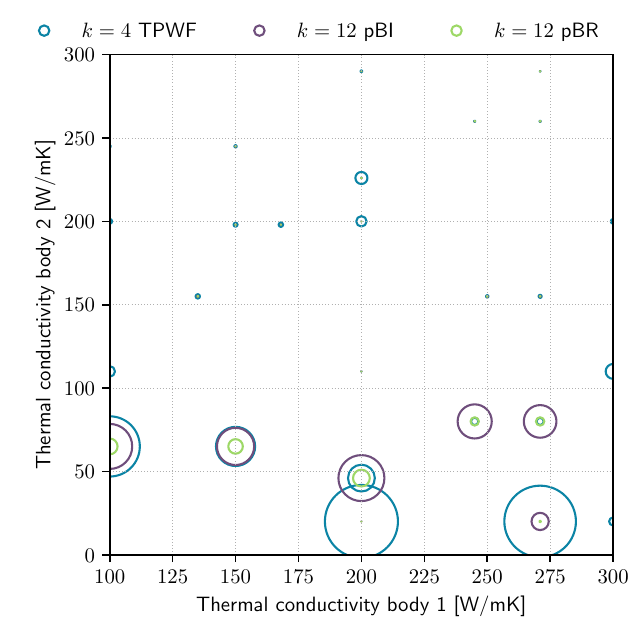}};
            \begin{scope}[x={(image.south east)},y={(image.north west)}]  
                \draw[red] (0.17,0.123) rectangle (0.952,0.42);
            \end{scope}
        \end{tikzpicture}
         \caption{Maximum MSRE for entire parameter space realized by the classical approach, pBI and pBR}
         \label{fig:PM_p2_Max_All}
     \end{subfigure}
     \hfill
     \begin{subfigure}[b]{0.46\textwidth}
         \centering
         \includegraphics[width=\textwidth]{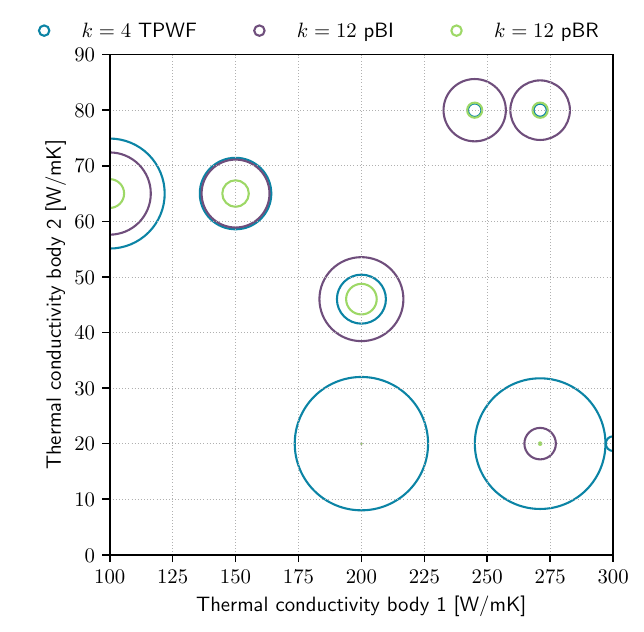}
         \caption{Maximum MSRE for $K_2=[20,110]\frac{\mathrm{W}}{\mathrm{mK}}$ realized by the classical approach, pBI and pBR}
         \label{fig:PM_p2_Max_All_Down}
     \end{subfigure}
        \caption{Bubble chart of the maximum mean square relative error of the power module two-parameter surrogate for the classical approach, parametric box interpolation and parametric box reduction. }
        \label{fig:PM_p2}
\end{figure}

Figure \ref{fig:PM_p2} presents the bubble charts of the pROMs temperature's mean square relative error for two regions. The x-coordinate denotes the thermal conductivity of body 1, $K_1$, and the y-coordinate is the thermal conductivity of body 2, $K_2$. The radius of each disk is proportional to the MSRE. Figure \ref{fig:PM_p2_Max_All} highlights the entire design space. The upper region with high values of both parameters shows smaller disks, corresponding to a low MSRE. Additionally, it can be observed that the error increases significantly for low values of $K_2$, this region is enclosed within a red rectangle. To investigate this behavior, Fig. \ref{fig:PM_p2_Max_All_Down} focuses on this region, i.e., $K_2 = [20,110] \frac{\mathrm{W}}{\mathrm{mK}}$. Contrary to the results for the beam model, which evidence a better performance with pBI and a higher number of training points, the pROM with $k = 12$  does not guarantee a more accurate model.

Therefore, the error arises during the basis change stage, which is not modified by pBI. As explained in Section \ref{ss:BR}, the surrogate model could lose accuracy when the system matrices' entries between training points have different orders of magnitude. For example, the following training points $\mu_{1} = [100,20]\frac{\mathrm{W}}{\mathrm{mK}}$ and $\mu_{9} = [300,290]\frac{\mathrm{W}}{\mathrm{mK}}$ have entries that differ by a factor of $14.51$; additionally, the matrices’ values of $\mu_{9}$ are $3.4$ times larger on average than the ones from $\mu_{1}$. During the basis change stage, the SVD artificially amplifies the influence of $\mu_{9}$ in the model because the most important modes of $\textbf{V}_{all}$ come from the training points with larger entries. Therefore, the transformation matrices and the global reduced space are affected. This drawback is amplified when more training points are added to the surrogate because a larger $\textbf{V}_{all}$ is constructed, but the $\textbf{R}$ matrix has the same dimension (Eq. \eqref{eq:RMatrix}). Thus, it reduces the influence of training models with low thermal conductivity values (e.g., $\mu_{1}$) even more. 

To remedy the issue, we applied pBR, as explained in Section \ref{ss:BR}. This is expected to reduce the non-physical influence of the outlying points on the surrogate model. To illustrate the accuracy of the novel technique, Fig. \ref{fig:PM_p2} depicts the temperature's maximum MSRE when pBR is employed. For $k=12$, the pBR improves the model's overall accuracy. Therefore, indicating that pBR can achieve higher accuracy when dealing with large parameter spaces and improve the approximation by increasing the number of training points.

\begin{figure}[ht]
     \centering
     \begin{subfigure}[b]{0.48\textwidth}
         \centering
         \includegraphics[width=\textwidth]{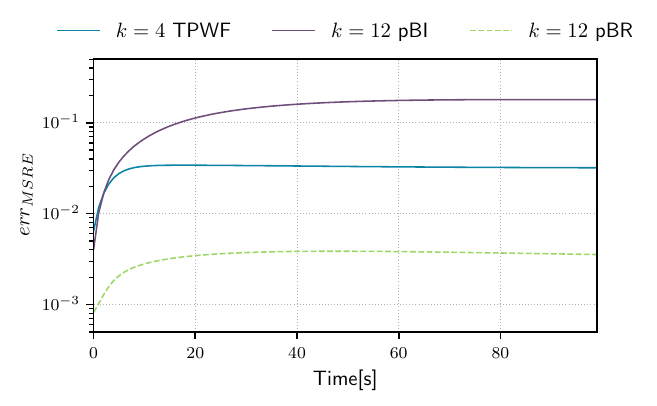}
         \caption{MSRE measured over the entire computational domain for $\mu_{13}$}
         \label{fig:PM_p2_MSRE_All_13}
     \end{subfigure}
     \hfill
     \begin{subfigure}[b]{0.48\textwidth}
         \centering
         \includegraphics[width=\textwidth]{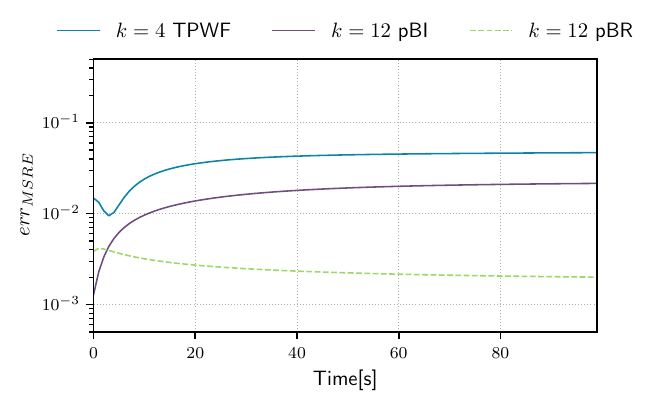}
         \caption{Local MSRE computed only on the chips temperature field for $\mu_{13}$}
         \label{fig:PM_p2_MSRE_Chip_13}
     \end{subfigure}
        \caption{Global and local mean square relative error analysis of pBI and pBR for $\mu_{13} = [271,80]\frac{\mathrm{W}}{\mathrm{mK}}$; pBR has the best overall accuracy and pBI can only accurately predict the chips' temperature}
        \label{fig:PM_p2_MSRE_13}
\end{figure}

In addition, two validation points are analyzed in more detail to better understand the behavior shown in Fig. \ref{fig:PM_p2}. For this purpose, we evaluate the local and global temperature predictions by investigating the temperature field MSRE for two domains. The first is the entire model's computational domain for a global analysis. Meanwhile, for the local analysis, the domain comprises only the chips (body 1 Fig. \ref{fig:PowerModule}), as the chips experience the decisive temperature gradients. Figure \ref{fig:PM_p2_MSRE_13} illustrates the temperature MSRE for $\mu_{13} = [271,80]\frac{\mathrm{W}}{\mathrm{mK}}$, a point with a relatively low $K_2$ and high $K_1$. The x-axis represents the number of time steps, whereas the y-axis measures the MSRE. The pBR is more accurate in both domains, having an MSRE of less than $0.5\%$ in the global and local approximation. When analyzing the global behavior depicted in Fig. \ref{fig:PM_p2_MSRE_All_13}, it is evident that sampling more training points without using pBR leads to higher errors. For $k = 4$ and $k = 12\; \mathrm{pBI}$, the global errors are $3.2\%$ and $17.9\%$, respectively. However, the behavior is the opposite when analyzing the error on the chips, as evidenced in Fig. \ref{fig:PM_p2_MSRE_Chip_13}. A surrogate model with $k = 12\; \mathrm{pBI}$ approximates the chips' temperature field more accurately than with only four training points.

To understand the results from Fig. \ref{fig:PM_p2_MSRE_13}, it is important to remember that both the Arnoldi method and the SVD try to identify the HFM dominant phenomena. For this application case, the dominant phenomena will be present in the chips because they have the highest temperature and temperature gradient. During the system reduction, the other power module components are inaccurately captured when the Krylov subspace is computed. Moreover, during the basis change stage, these components are again disregarded for training points whose matrix entries are significantly smaller than the other points. In addition, when increasing $k$, their influence decreases even more, as previously explained. Therefore, the pROM for $ k=12\; \mathrm{pBI}$ has a greater MSRE over the entire computational domain for $\mu_{13}$. In contrast, the behavior of high-temperature regions (i.e., chips) is better conserved after the system reduction and basis change. This allows a more accurate approximation of the chips' temperature field with more training points, even when using $k=12\; \mathrm{pBI}$. For these reasons, pROMs generated with pBR have a lower global and local error.

\begin{figure}[ht]
     \centering
     \begin{subfigure}[b]{0.48\textwidth}
         \centering
         \includegraphics[width=\textwidth]{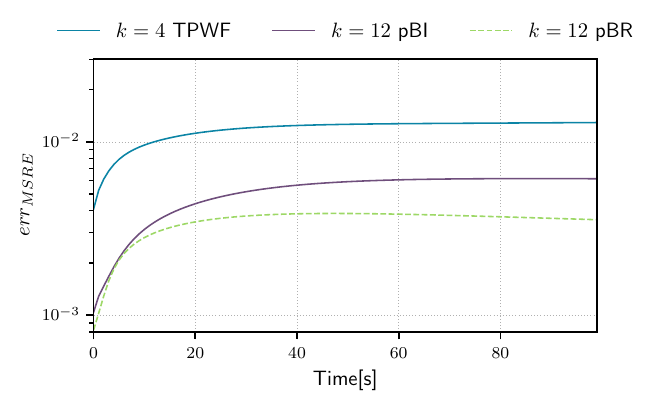}
         \caption{MSRE measured over the entire computational domain for $\mu_{24}$}
         \label{fig:PM_p2_MSRE_All_24}
     \end{subfigure}
     \hfill
     \begin{subfigure}[b]{0.48\textwidth}
         \centering
         \includegraphics[width=\textwidth]{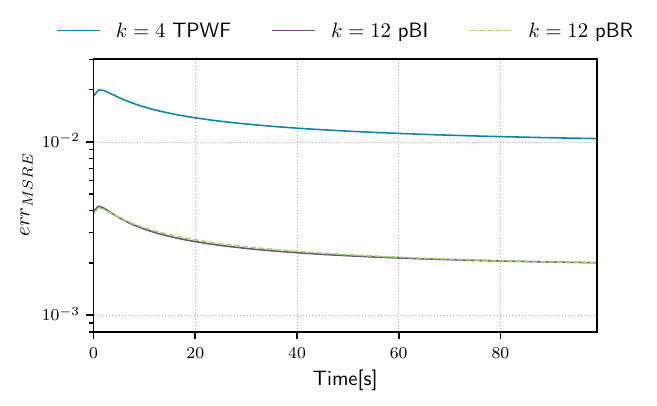}
         \caption{Local MSRE computed only on the chips temperature field for $\mu_{24}$}
         \label{fig:PM_p2_MSRE_Chip_24}
     \end{subfigure}
        \caption{Global and local mean square relative error analysis of pBI and pBR for $\mu_{24} = [150,245]\frac{\mathrm{W}}{\mathrm{mK}}$. Both pBI and pBR have a improve accuracy with more training points in both global and local predictions}
        \label{fig:PM_p2_MSRE_24}
\end{figure}

Similarly, Fig. \ref{fig:PM_p2_MSRE_24} evaluates the MSRE for $\mu_{24} = [150,245]\frac{\mathrm{W}}{\mathrm{mK}}$. The x-axis represents the number of time steps, whereas the y-axis measures the MSRE. For the previous case, the training points used for the matrix interpolation have the lowest influence on the basis change; this is not the case for $\mu_{24}$ because both $K_1$ and $K_2$ have a relatively high value, similar to a set up in Fig. \ref{fig:pBI_k9H}. Thus, the training points with the longest distance to $\mu_{24}$ have a negligible negative impact on the global approximation than that evidenced for $\mu_{13}$. To illustrate this, Fig. \ref{fig:PM_p2_MSRE_All_24} shows the MSRE over the entire computational domain. Without pBR, $k = 12$ has a lower error than $k = 4$. By applying pBR, we can increase the global accuracy by $42\%$. For the local analysis, the MSRE of the chips is depicted in Fig. \ref{fig:PM_p2_MSRE_Chip_24}. As with the global error, the surrogate with four training points has the highest MSRE. The pBR improves only negligibly due to the low non-physical influence of the other training points and components in the chips' temperature prediction. 

\begin{figure}[!ht]
    \centering
    \begin{tikzpicture}
        \node[anchor=south west,inner sep=0] (image) at (0,0) {\includegraphics[width=0.5\linewidth]{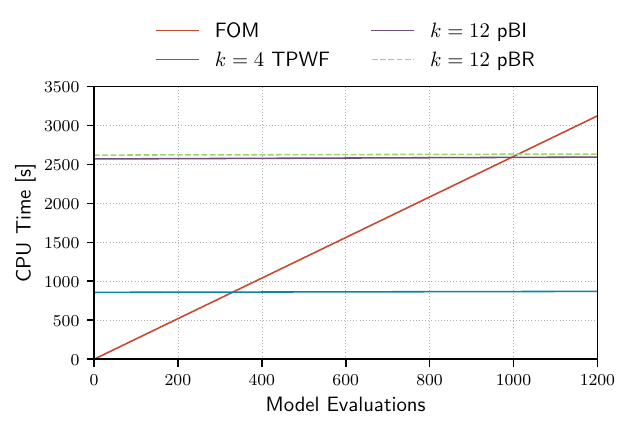}};
        \begin{scope}[x={(image.south east)},y={(image.north west)}]  
            \tikzstyle{arrow} = [->,>=stealth]
            \draw[arrow] (0.8,0.5) -- (0.81,0.63);
            \node[black] at (0.8,0.45) {\shortstack{\tiny\textbf{Break-even} \\ \tiny \textbf{point}}};
            \draw[arrow] (0.7,0.45) -- (0.39,0.33);
        \end{scope}
    \end{tikzpicture}
    \caption{Computational time comparison between the FOM and the different surrogate models of the power module. The increase of computational time between pBI and pBR is negligible}
    \label{fig:PM_time}
\end{figure}

As before, we evaluate the computational cost of the proposed technique. Figure \ref{fig:PM_time} shows the computational effort needed to create the two-parameter surrogate model. The horizontal axis illustrates the number of model evaluations. For all HFM, the Krylov subspace method to create the local ROMs is the same. Therefore, the required time for the system reduction depends only on the number of training points. For the classical approach and pBI, only one SVD is computed for the matrix transformation.

Meanwhile, the pBR requires multiple SVD computations. These can be either done during the online or offline phase, as described in Section \ref{ss:BR}. For this work, the analysis is carried out by computing the SVDs during the offline stage. Six SVDs are carried out, one for each box in the $3\times 4$ training grid. Although the five additional SVDs, the basis change time for $k=12 \; \mathrm{pBR}$ only increases $23\%$ compared to $k=12 \; \mathrm{pBI}$. This behavior is due to the additional SVD being executed on matrices of smaller dimensions, and the SVD complexity grows nonlinearly \cite{li2019tutorialcomplexityanalysissingular}. Furthermore, this additional computational effort is minimal, constituting just a $2\%$ increment in the offline time because the predominant computational cost is associated with the system reduction, which does not change for the same number of training points as explained above. This difference is highlighted in the y-axis intercept of the plot depicted in Fig. \ref{fig:PM_time}. Similar to the beam, the speed-up between the methods remains effectively the same. However, the speed-up for the power module is $225$, 15 times faster than for the beam because the original system has a higher number of DoFs. After only $1011$ time integration steps, the pBR requires less time than the FOM.
\section{Conclusion}\label{ss:conclusion}

In this contribution, we introduce a novel matrix interpolation technique for parametric reduced order models. Specifically, we propose two new algorithms: box clustering and tensor product weight function. A noteworthy characteristic is that these methods remove the parameter space normalization step, which, in many cases, affects the surrogate accuracy. Moreover, box clustering overcomes the challenges of the nearest-neighbor strategy for high-dimensional spaces. By integrating these two new schemes, we introduce pBI and pBR, two reduction methods that minimize the non-physical numerical artifacts when handling large parameter spaces. The methods are divided into three main stages: (i) system reduction, where with a set of training points reduced models are created with a Krylov Subspace Method; (ii) basis change, where all reduced models are projected into a global reduced space; (iii) matrix interpolation, where the new reduced matrices for any given combination of input parameters are computed with the tensor product weight function. Instead of using all training points, the box clustering creates a subset of these points that are used for (ii) and (iii). 

Notably, the experiments demonstrate that increasing the number of training points does not guarantee a lower error for the classical approach. Despite pBI demonstrating promising results for the 3D beam, it faces challenges when the dimension of the FOM increases for the power module. For such cases, the numerical experiments demonstrate that pBI can accurately predict the high-temperature regions but fails to approximate the system's global behavior. Whereas, pBR accurately predicts the global and local temperature field for the power module and the 3D beam. 

Furthermore, the additional computational effort required for the new technique is negligible. Nevertheless, one of the main drawbacks of the method is that it suffers from the curse of dimensionality. A formulation for higher dimensions with unstructured training grids is still an open question. To conclude, through an academic and industrial application case, we establish that pBR can successfully create surrogate models from a wide range of material properties for structural and thermal problems. The prediction achieves higher accuracy than the classical approach without significantly impacting the CPU time. In future work, we will explore issues with varying geometrical parameters, specifically when the mesh changes between parameter points.

\section{CRediT authorship contribution statement}
\textbf{Juan Angelo Vargas-Fajardo:} Conceptualization, Data curation, Formal analysis, Investigation, Methodology, Software, Validation, Visualization, Writing  – original draft.
\textbf{Diana Manvelyan:} Conceptualization, Methodology, Supervision, Writing – review and editing, Project administration, Resources.
\textbf{Catharina Czech:} Conceptualization, Methodology, Supervision, Writing – review and editing, Project administration, Resources.
\textbf{Pietro Botazzoli:} Conceptualization, Supervision, Writing – review and editing, Resources.
\textbf{Fabian Duddeck:} Supervision, Writing – review and editing.

\section{Funding sources}
The work was supported by Siemens AG.

\section{Declaration of competing interest}
The authors declare that they have no known competing financial interests or personal relationships that could have appeared to influence the work reported in this paper.

\section{Acknowledgements}
The authors would like to thank Christoph Heinrich and his research group for valuable discussions

\bibliographystyle{unsrt} 
\bibliography{references}

\end{document}